\documentclass[reqno,12pt,a4paper]{amsart}
\usepackage{amssymb}
\def\C{\mathbb{C}}
\def\R{\mathbb{R}}

\def\CP{\mathbb{CP}}
\def\dbar{\bar{\partial}}
\newtheorem{thm}{Theorem}[section]

\newtheorem{prop}[thm]{Proposition}

\theoremstyle{definition}
\newtheorem{defn}{Definition}[section]
\theoremstyle{remark}
\newtheorem{rem}{\rm\bfseries{Remark}}[section]

\title{Symplectic maps to projective spaces and symplectic invariants}
\author{Denis Auroux}
\address{Centre de Math\'ematiques, Ecole Polytechnique, F-91128 Palaiseau,
France}
\email{auroux@math.polytechnique.fr}
\begin{document}

\begin{abstract}
After reviewing recent results on symplectic Lefschetz pencils and
symplectic branched covers of $\CP^2$, we describe a new construction of
maps from symplectic manifolds of any dimension to $\CP^2$ and the
associated monodromy invariants. We also show that a dimensional induction
process makes it possible to describe any compact symplectic manifold by
a series of words in braid groups and a word in a symmetric group.
\end{abstract}

\maketitle

\section{Introduction}

Let $(X^{2n},\omega)$ be a compact symplectic manifold. We will throughout
this text assume that the cohomology class $\frac{1}{2\pi}[\omega]\in
H^2(X,\R)$ is integral. This assumption makes it possible to define a
complex line bundle $L$ over $X$ such that $c_1(L)=\frac{1}{2\pi}[\omega]$.
We also endow $X$ with a compatible almost-complex structure $J$, and endow
$L$ with a Hermitian metric and a Hermitian connection of curvature
$-i\omega$.

The line bundle $L$ should be thought of as a symplectic version of an ample
line bundle over a complex manifold. Indeed, although the lack of
integrability of $J$ prevents the existence of holomorphic sections, it
was observed by Donaldson in \cite{D1} that, for large $k$, the line 
bundles $L^{\otimes k}$ admit many approximately holomorphic sections.

Observe that all results actually apply as well to the case 
where $\frac{1}{2\pi}[\omega]$ is not integral, with the only difference that
the choice of the line bundle $L$ is less natural~: the idea is to perturb
$\omega$ into a symplectic form $\omega'$ whose cohomology class is
rational, and then work with a suitable multiple of $\omega'$. One chooses
an almost-complex structure $J'$ which simultaneously is compatible
with $\omega'$ and satisfies the positivity property $\omega(v,J'v)>0$ for
all tangent vectors. All the objects that we construct are then
approximately $J'$-holomorphic, and therefore symplectic with respect to
not only $\omega'$ but also $\omega$.

Donaldson was the first to show in \cite{D1} that, among the many
approximately holomorphic sections of $L^{\otimes k}$ for $k\gg 0$, there
is enough flexibility in order to obtain nice transversality properties~;
this makes it possible to imitate various classical topological
constructions from complex algebraic geometry in the symplectic category.
Let us mention in particular the construction of smooth symplectic
submanifolds (\cite{D1}, see also \cite{A1} and \cite{presas}), 
symplectic Lefschetz pencils (\cite{D2}, see also \cite{Dsurvey}),
branched covering maps to $\CP^2$ (\cite{A2},\cite{AK}), Grassmannian
embeddings and determinantal submanifolds (\cite{presas}).
\medskip

Intuitively, the main reason why the approximately holomorphic framework
is suitable to imitate results from algebraic geometry is that, for large
values of $k$, the increasing curvature of $L^{\otimes k}$ provides access
to the geometry of $X$ at very small scale~; as one zooms into $X$, the
geometry becomes closer and closer to a standard complex model, and the
lack of integrability of $J$ becomes negligible.

In all cases the strategy is more or less the same~: the aim is, starting
with a sequence of given approximately holomorphic sections of $L^{\otimes k}$
for all $k\gg 0$, to perturb them in order to ensure uniform transversality
properties that will guarantee the desired topological properties. 

For example, to construct submanifolds, one obtains bounds of the type
$|\nabla s_k|_{g_k}>\eta$ along the zero set of $s_k$ for a fixed constant
$\eta>0$ independent of $k$, while approximate holomorphicity implies a
bound of the type $|\dbar s_k|_{g_k}=O(k^{-1/2})$ everywhere. Here $g_k=kg$
is a rescaled metric which dilates everything by a factor of $k^{1/2}$ in
order to adapt to the decreasing ``characteristic scale'' imposed by the
increasing curvature $-ik\omega$ of the line bundles $L^{\otimes k}$.
The desired topological picture, similar to the complex algebraic case, 
emerges for large $k$ as an inequality of the form $|\dbar s_k|\ll 
|\partial s_k|$ becomes satisfied~: this can easily be shown to imply that 
the zero set of $s_k$ is smooth, approximately pseudo-holomorphic, and
symplectic.

The starting points for the construction, in all cases, are
the existence of very localized approximately holomorphic sections of
$L^{\otimes k}$ concentrated near any given point $x\in X$, and an effective
transversality result for approximately holomorphic functions defined over
a ball in $\C^n$ with values in $\C^r$ due to Donaldson (see \cite{D1} for 
the case $r=1$ and \cite{D2} for the general case). These two ingredients
imply that a small localized perturbation can be used to ensure uniform
transversality over a small ball. Combining this local result with
a globalization argument (\cite{D1}, see also \cite{A2} and \cite{presas}),
one obtains transversality everywhere.
\medskip

The interpretation of the construction of submanifolds as an
effective transversality result for sections extends verbatim to the more 
sophisticated constructions (Lefschetz pencils, branched coverings)~:
in these cases the transversality properties also concern the covariant
derivatives of the sections, and this can be thought of as an effective
analogue in the approximately holomorphic category of the standard 
generalized transversality theorem for jets.

This is especially clear when looking at the arguments in \cite{presas},
\cite{D2} or \cite{A2}~: the perturbative argument is now used
to obtain uniform transversality of the holomorphic parts of the 
$1$-jets or $2$-jets of the sections with respect to certain closed 
submanifolds in the space of holomorphic jets.
Successive perturbations are used to obtain transversality to the various
strata describing the possible singular models~; one uses that each stratum
is smooth away from lower dimensional strata, and that transversality to
these lower dimensional strata is enough to imply transversality to
the higher dimensional stratum near its singularities.

An extra step is necessary in the constructions~: recall that desired
topological properties only hold when the antiholomorphic parts of the
derivatives are much smaller than the holomorphic parts. In spite of
approximate holomorphicity, this can be a problem when the holomorphic part
of the jet becomes singular. Therefore, a small perturbation is needed to 
kill the antiholomorphic part of the jet near the singularities~; this
perturbation is in practice easy to construct. The reader is referred to
\cite{D2} and \cite{A2} for details.

Although no general statement has yet been formulated and proved, it is
completely clear that a very general result of uniform 
transversality for jets holds in the approximately holomorphic category.
Therefore, the observed phenomenon for Lefschetz pencils and maps to
$\CP^2$, namely the fact that near every point $x\in X$ the constructed
maps are given in approximately holomorphic coordinates by one of the
standard local models for generic holomorphic maps, should hold in all
generality, independently of the dimensions of the source and target spaces.
This approach will be developed in a forthcoming paper \cite{Atrans}.
\medskip

In the remainder of this paper we focus on the topological monodromy
invariants that can be derived from the various available constructions.
In Section 2 we study symplectic Lefschetz pencils and their monodromy,
following the results of Donaldson \cite{D2} and Seidel \cite{seidel}. 
In Section 3 we describe symplectic branched covers of $\CP^2$ and their 
monodromy invariants, following \cite{A2} and \cite{AK}~; we also discuss the
connection with 4-dimensional Lefschetz pencils. In Section 4 we extend this
framework to the higher dimensional case, and investigate a new type of
monodromy invariants arising from symplectic maps to $\CP^2$. We finally
show in Section 5 that a dimensional induction process makes it possible
to describe a compact symplectic manifold of any dimension by
a series of words in braid groups and a word in a symmetric group.
\medskip

{\bf Acknowledgement.} The author wishes to thank Ludmil Katzarkov,
Paul Seidel and Bob Gompf for stimulating discussions, as well as 
Simon Donaldson for his interest in this work.

\section{Symplectic Lefschetz pencils}

Let $(X^{2n},\omega)$ be a compact symplectic manifold as above, and let
$s_0,s_1$ be suitably chosen approximately
holomorphic sections of $L^{\otimes k}$. Then $X$ is endowed with a
structure of {\em symplectic Lefschetz pencil}, which can be described as 
follows.

For any $\alpha\in\CP^1=\C\cup\{\infty\}$, define 
$\Sigma_\alpha=\{x\in X,\ s_0+\alpha s_1=0\}$. Then the submanifolds
$\Sigma_\alpha$ are symplectic hypersurfaces, smooth except for finitely 
many values of the parameter $\alpha$~; for these parameter values 
$\Sigma_\alpha$ contains a singular point (a normal crossing when $\dim
X=4$). Moreover, the submanifolds $\Sigma_\alpha$ fill all of $X$,
and they intersect transversely along a codimension 4 symplectic 
submanifold $Z=\{x\in X,\ s_0=s_1=0\}$, called the set of {\em base points} 
of the pencil.

Define the projective map $f=(s_0\!:\!s_1):X-Z\to\CP^1$, whose level sets
are precisely the hypersurfaces $\Sigma_\alpha$. Then $f$ is required to be 
a {\em complex Morse function}, i.e.\ its critical points are
isolated and non-degenerate, with local model $f(z_1,\dots,z_n)=z_1^2+\dots+
z_n^2$ in approximately holomorphic coordinates.

The following result due to Donaldson holds~:

\begin{thm}[Donaldson \cite{D2}]
For $k\gg 0$, two suitably chosen approximately holomorphic sections of 
$L^{\otimes k}$ endow $X$ with a structure of symplectic Lefschetz pencil, 
canonical up to isotopy.
\end{thm}

This result is proved by obtaining uniform transversality with respect to 
the strata $s_0=s_1=0$ (of complex codimension 2) and $\partial f=0$ (of 
complex codimension $n$) in the space of holomorphic 1-jets of sections of
$\C^2\otimes L^{\otimes k}$, by means of the techniques described in the
introduction. A small additional perturbation ensures the compatibility
requirement that $\dbar f$ vanishes at the points where $\partial f=0$. 
These properties are sufficient to ensure that the structure is that of a 
symplectic Lefschetz pencil. For details, the reader is referred to \cite{D2}.

The statement that the constructed pencils are canonical up to isotopy for
$k\gg 0$ is to be interpreted as follows. Consider two sequences
$(s^0_k)_{k\gg 0}$ and $(s^1_k)_{k\gg 0}$ of 
approximately holomorphic sections of $\C^2\otimes L^{\otimes k}$ for
increasing values of $k$. Assume that they satisfy the three above-described
transversality and compatibility properties and hence define symplectic 
Lefschetz pencils. Then, for large enough $k$ (how large exactly depends 
on the estimates on the given sections), there exists an interpolating family
$(s^t_k)_{t\in [0,1]}$ of approximately holomorphic sections, depending
continuously on the parameter $t$, such that for all values of $t$ the
sections $s^t_k$ satisfy the transversality and compatibility properties. 
In particular, for large enough $k$ the symplectic Lefschetz pencils
defined by $s^0_k$ and $s^1_k$ are isotopic to each other. Moreover, the
same result remains true if the almost-complex structures $J_0$ and $J_1$
with respect to which $s^0_k$ and $s^1_k$ are approximately holomorphic
differ, so the topology of the constructed pencils depends only on the
topology of the symplectic manifold $X$ (and on $k$ of course). However, 
because isotopy holds only for large values of $k$, this is only a weak 
(asymptotic) uniqueness result.
\medskip

A convenient way to study the topology of a Lefschetz pencil is to blow up
$X$ along the submanifold $Z$. The resulting symplectic manifold $\hat{X}$
is the total space of a {\em symplectic Lefschetz fibration} 
$\hat{f}:\hat{X}\to\CP^1$.  Although in the following description
we work on the blown up manifold $\hat{X}$, it is actually preferrable to
work directly on $X$~; verifying that the discussion applies to $X$ itself
is a simple task left to the reader.

The fibers of $\hat{f}$ can be identified with
the submanifolds $\Sigma_\alpha$, made mutually disjoint by the blow-up
process. It is then possible to study the {\em monodromy} of the fibration
$\hat{f}$ around its singular fibers.

One easily checks that this monodromy
consists of symplectic automorphisms of the fiber $\Sigma_\alpha$.
Moreover, the exceptional divisor obtained by blowing up the set of 
base points $Z$ is a subfibration of $\hat{f}$, with fiber $Z$, which
is unaffected by the monodromy. Therefore, it is natural to consider that
the monodromy of $\hat{f}$ takes values in the symplectic mapping class
group $\mathrm{Map}^\omega(\Sigma,Z)=\pi_0(\{\phi\in
\mathrm{Symp}(\Sigma,\omega),\phi_{|U(Z)}=\mathrm{Id}\})$, i.e.\ the
set of isotopy classes of symplectomorphisms of the generic fiber $\Sigma$
which coincide with the identity near $Z$. 

In the four-dimensional case, $Z$ consists of a finite number $n$ of points,
and $\Sigma$ is a compact surface with a certain genus $g$ (note that
$\Sigma$ is always connected because it satisfies a Lefschetz hyperplane
type property)~; $\mathrm{Map}^\omega(\Sigma,Z)$ is then the classical
mapping class group $\mathrm{Map}_{g,n}$ of a genus $g$ surface with $n$
boundary components.

In fact, the image of the monodromy map is contained in the subgroup of
{\it exact} symplectomorphisms in $\mathrm{Map}^\omega(\Sigma,Z)$~: the
connection on $L^{\otimes k}$ induces over $\Sigma-Z$ a 1-form $\alpha$
such that $d\alpha=\omega$. This endows $\Sigma-Z$ with a structure of 
exact symplectic manifold. Monodromy transformations are then exact 
symplectomorphisms in the sense that they preserve not only $\omega$ but
also the 1-form $\alpha$ (see \cite{seidel2} for details).
\medskip

It is well-known (see e.g.\ \cite{seidel}, \cite{seidel2}) that the singular 
fibers of a Lefschetz fibration are obtained from the generic fiber by 
collapsing a {\em vanishing cycle} to a point. The vanishing cycle is an 
embedded closed loop in $\Sigma$ in the four-dimensional case~; more 
generally, it is an embedded Lagrangian sphere $S^{n-1}\subset \Sigma$.
Then, the monodromy of $\hat{f}$ around one of its singular fibers consists 
in a {\em generalized Dehn twist} in the positive direction along the 
vanishing cycle.

The picture is the following~:

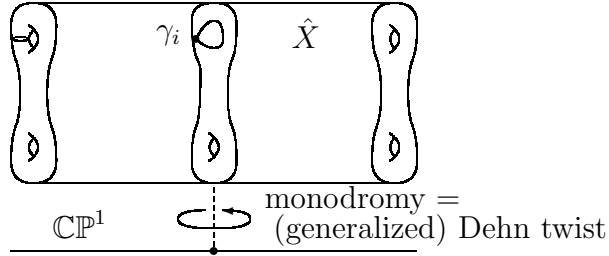
\setlength{\unitlength}{0.6mm}
\begin{center}
\begin{picture}(80,65)(-40,-20)

\put(9,-5.5){ monodromy =}
\put(11,-12){ (generalized) Dehn twist}
\qbezier(-2,-6)(-8,-6)(-8,-8)
\qbezier(-8,-8)(-8,-10)(0,-10)
\qbezier(0,-10)(8,-10)(8,-8)
\qbezier(8,-8)(8,-6)(2,-6)
\put(2,-6){\vector(-1,0){0}}

\put(-40,0){\line(1,0){80}}
\put(-40,40){\line(1,0){80}}
\put(-45,-15){\line(1,0){90}}
\put(-38,-12){ $\CP^1$}
\qbezier(-40,40)(-35,40)(-35,35)
\qbezier(-35,35)(-35,30)(-36,27.5)
\qbezier(-36,27.5)(-37,25)(-37,20)
\qbezier(-40,0)(-35,0)(-35,5)
\qbezier(-35,5)(-35,10)(-36,12.5)
\qbezier(-36,12.5)(-37,15)(-37,20)
\qbezier(-40,40)(-45,40)(-45,35)
\qbezier(-45,35)(-45,30)(-44,27.5)
\qbezier(-44,27.5)(-43,25)(-43,20)
\qbezier(-40,0)(-45,0)(-45,5)
\qbezier(-45,5)(-45,10)(-44,12.5)
\qbezier(-44,12.5)(-43,15)(-43,20)
\qbezier(-41,35)(-37,32)(-41,29)
\qbezier(-40,34)(-42,32)(-40,30)
\qbezier(-41,5)(-37,8)(-41,11)
\qbezier(-40,6)(-42,8)(-40,10)

\qbezier(40,40)(35,40)(35,35)
\qbezier(35,35)(35,30)(36,27.5)
\qbezier(36,27.5)(37,25)(37,20)
\qbezier(40,0)(35,0)(35,5)
\qbezier(35,5)(35,10)(36,12.5)
\qbezier(36,12.5)(37,15)(37,20)
\qbezier(40,40)(45,40)(45,35)
\qbezier(45,35)(45,30)(44,27.5)
\qbezier(44,27.5)(43,25)(43,20)
\qbezier(40,0)(45,0)(45,5)
\qbezier(45,5)(45,10)(44,12.5)
\qbezier(44,12.5)(43,15)(43,20)
\qbezier(39,35)(43,32)(39,29)
\qbezier(40,34)(38,32)(40,30)
\qbezier(39,5)(43,8)(39,11)
\qbezier(40,6)(38,8)(40,10)

\qbezier(0,40)(5,40)(5,35)
\qbezier(5,35)(5,30)(4,27.5)
\qbezier(4,27.5)(3,25)(3,20)
\qbezier(0,0)(5,0)(5,5)
\qbezier(5,5)(5,10)(4,12.5)
\qbezier(4,12.5)(3,15)(3,20)
\qbezier(0,40)(-5,40)(-5,36)
\qbezier(-5,36)(-5,32)(-4,28)
\qbezier(-4,28)(-3,26)(-3,25)
\qbezier(0,0)(-5,0)(-5,5)
\qbezier(-5,5)(-5,10)(-4,12.5)
\qbezier(-4,12.5)(-3,15)(-3,25)
\qbezier(2,32)(2,36)(0,36)
\qbezier(0,36)(-2,36)(-4,32)
\qbezier(2,32)(2,30)(0,30)
\qbezier(0,30)(-2,30)(-4,32)
\qbezier(-1,5)(3,8)(-1,11)
\qbezier(0,6)(-2,8)(0,10)

\put(15,30){ $\hat{X}$}

\put(-4,32){\circle*{2}}
\put(-15,32){ $\gamma_i$}
\qbezier(-44.8,32)(-42.9,31)(-41,32)
\qbezier[8](-44.8,32)(-42.9,33)(-41,32)
\multiput(0,-1)(0,-2){7}{\line(0,-1){1}}
\put(0,-15){\circle*{2}}

\end{picture}
\end{center}

However, because the normal bundle to the exceptional divisor is
not trivial, the monodromy map cannot be defined over all of $\CP^1$, and
we need to restrict ourselves to the preimage of an affine subset $\C$
(the fiber at infinity can be assumed regular). The monodromy around the
fiber at infinity of $\hat{f}$ is given by a mapping class group element
$\delta_Z$ corresponding to a {\em twist around $Z$}. In the
four-dimensional case $Z$ consists of $n$ points, and $\delta_Z$ is the
product of positive Dehn twists along $n$ loops each encircling one of
the base points~; in the higher-dimensional case $\delta_Z$ is a positive 
Dehn twist along the unit sphere bundle in the normal bundle of $Z$ in 
$\Sigma$ (i.e.\ it restricts to each fiber of the normal bundle as a Dehn
twist around the origin).

It follows from the above observations 
that the monodromy of the Lefschetz fibration $\hat{f}$
with critical levels $p_1,\dots,p_d$ is given by a group homomorphism
\begin{equation}
\psi:\pi_1(\C-\{p_1,\dots,p_d\})\to\mathrm{Map}^\omega(\Sigma^{2n-2},Z)
\end{equation}
which maps the {\em geometric generators} of $\pi_1(\C-\{p_1,\dots,p_d\})$,
i.e.\ loops going around one of the points $p_i$, to Dehn twists.

Alternately, choosing a system of generating loops in
$\C-\{p_1,\dots,p_d\}$, we can express the monodromy by a {\em factorization}
of $\delta_Z$ in the mapping class group~:
\begin{equation}
\delta_Z=\prod_{i=1}^d \tau_{\gamma_i},\end{equation} 
where $\gamma_i$ is the image
in a chosen reference fiber of the vanishing cycle of the singular fiber 
above $p_i$ and $\tau_{\gamma_i}$ is the corresponding positive Dehn twist. 
The identity (2) in $\mathrm{Map}^\omega(\Sigma,Z)$ expresses the fact that 
the monodromy of the fibration around the point at infinity in $\CP^1$ 
decomposes as the product of the elementary monodromies around each of 
the singular fibers.

The monodromy morphism (1), or equivalently the mapping class group 
factorization (2), completely characterize the topology of the Lefschetz
fibration $\hat{X}$. However, they are not entirely canonical, because
two choices have been implicitly made in order to define them. 

First, a base point in $\C-\{p_1,\dots,p_d\}$ and an identification 
symplectomorphism between $\Sigma$ and the chosen reference fiber of 
$\hat{f}$ are needed in order to view the monodromy transformations 
as elements in the mapping class group of $\Sigma$.  The choice of a 
different identification affects the monodromy morphism
$\psi$ by conjugation by a certain element 
$g\in\mathrm{Map}^\omega(\Sigma,Z)$. The corresponding operation on the
mapping class group factorization (2) is a {\it simultaneous conjugation}
of all factors~: each factor $\tau_{\gamma_i}$ is replaced by
$\tau_{g(\gamma_i)}=g^{-1}\tau_{\gamma_i}g$.

Secondly, a system of generating loops has to be chosen in order to define 
a factorization of $\delta_Z$. Different choices of generatic systems differ
by a sequence of {\it Hurwitz operations}, i.e.\ moves in which two
consecutive generating loops are exchanged, one of them being conjugated by
the other in order to preserve the counterclockwise ordering. On the level
of the factorization, this amounts to replacing two consecutive factors
$\tau_1$ and $\tau_2$ by respectively $\tau_2$ and $\tau_2^{-1}\tau_1\tau_2$
(or, by the reverse operation, $\tau_1\tau_2\tau_1^{-1}$ and $\tau_1$).

It is quite easy to see that any two factorizations of $\delta_Z$ describing
the Lefschetz fibration $\hat{f}$ differ by a sequence of these two 
operations (simultaneous conjugation and Hurwitz moves). Therefore,
Donaldson's uniqueness statement implies that, for large enough values of
$k$, the mapping class group factorizations associated to the symplectic
Lefschetz pencil structures obtained in Theorem 2.1 are, up to simultaneous
conjugation and Hurwitz moves, {\it symplectic invariants} of the manifold
$(X,\omega)$.

Conversely, given any factorization of $\delta_Z$ in 
$\mathrm{Map}^\omega(\Sigma,Z)$ as a product of positive Dehn twists, it
is possible to construct a symplectic Lefschetz fibration with the given
monodromy. It follows from a result of Gompf that the total space of such
a fibration is always a symplectic manifold. In fact, because the monodromy
preserves the symplectic submanifold $Z\subset\Sigma$, it is also possible
to reconstruct the blown down manifold $X$. More precisely, the following
result holds~:

\begin{thm}[Gompf]
Let $(\Sigma,\omega_\Sigma)$ be a compact symplectic manifold, and 
$Z\subset\Sigma$ a codimension $2$ symplectic submanifold such that 
$[Z]=PD([\omega_\Sigma])$. Consider a factorization of $\delta_Z$ as a
product of positive Dehn twists in $\mathrm{Map}^\omega(\Sigma,Z)$.
In the case $\dim(\Sigma)=2$, assume moreover that all the Dehn twists
in the factorization are along loops that are not homologically trivial in
$\Sigma-Z$.

Then the total space $X$ of the corresponding Lefschetz pencil carries
a symplectic form $\omega_X$ such that, given a generic fiber $\Sigma_0$
of the pencil, $[\omega_X]$ is Poincar\'e dual to $[\Sigma_0]$, and
$(\Sigma_0,\omega_{X|\Sigma_0})$ is symplectomorphic to
$(\Sigma,\omega_\Sigma)$.
This symplectic structure on $X$ is canonical up to symplectic isotopy.
\end{thm}

The strategy of proof is to first construct a symplectic structure in the
correct cohomology class on a neighborhood of any fiber of the pencil, which
is easily done as $\Sigma$ already carries a symplectic structure and the
monodromy lies in the symplectomorphism group. It is then possible to
combine these local symplectic structures and obtain a globally defined
symplectic form, singular near the base locus $Z$. Since the total
monodromy is $\delta_Z$, the structure of $X$ near $Z$ is
completely standard, and so a non-singular symplectic form on $X$ can be 
recovered (this process can also be viewed as a symplectic blow-down along
the exceptional hypersurface $\CP^1\times Z$ in the total space of the
corresponding Lefschetz fibration). This operation changes the cohomology
class of the symplectic form on $X$, but one easily checks that the
resulting class is a nonzero multiple of the Poincar\'e dual to a fiber~;
scaling the symplectic form by a suitable factor then yields $\omega_X$.
The proof that this process is canonical up to symplectic isotopy is a
direct application of Moser's stability theorem. The reader is referred to
\cite{hypergompf} for details.
\medskip


In conclusion, the study of the monodromy of symplectic Lefschetz pencils 
makes it possible to define invariants of compact symplectic
manifolds, which in principle provide a complete description of the
topology. However, the complexity of mapping class groups and the
difficulties in computing the invariants in concrete situations greatly
decrease their usefulness in practice. This motivates the introduction of
other similar topological constructions which may lead to more usable
invariants.

\section{Branched covers of $\CP^2$ and invariants of symplectic
4-manifolds}

Throughout \S 3, we assume that $(X,\omega)$ is a compact symplectic
4-manifold. In that case, three generic approximately holomorphic 
sections $s_0$, $s_1$ and $s_2$ of $L^{\otimes k}$ never vanish
simultaneously, and so they define a projective map
$f=(s_0\!:\!s_1\!:\!s_2):X\to\CP^2$. It was shown in \cite{A2} that, if
the sections are suitably chosen, this map is a {\it branched covering},
whose branch curve $R\subset X$ is a smooth connected symplectic submanifold
in $X$. 

There are two possible local models in approximately holomorphic
coordinates for the map $f$ near the branch curve. The first one,
corresponding to a generic point of $R$, is the map 
$(x,y)\mapsto (x^2,y)$~; locally, both the branch curve $R$ and its image by
$f$ are smooth. The other local
model corresponds to the isolated points where $f$ does not restrict
to $R$ as an immersion. The model map is then $(x,y)\mapsto (x^3-xy,y)$,
and the image of the smooth branch curve $R:3x^2-y=0$ has equation
$f(R):27z_1^2=4z_2^3$ and presents a cusp singularity. These two local
models are the same as in the complex algebraic setting.

It is easy to see by considering the two model maps that $R$ is a smooth 
approximately holomorphic (and therefore symplectic) curve in $X$, and that 
$f(R)$ is an approximately holomorphic symplectic curve in $\CP^2$, immersed 
away from its cusps. After a generic perturbation, we can moreover require
that the branch curve $D=f(R)$ satisfies a self-transversality property,
i.e.\ that its only singular points besides the cusps are transverse double
points (``nodes''). Even though $D$ is approximately holomorphic, it is not 
immediately possible to require that all of its double points correspond 
to a positive intersection number with respect to the standard orientation 
of $\CP^2$~; the presence of (necessarily badly transverse) negative double 
points is a priori possible.

It was also shown in \cite{A2} that the branched coverings obtained from
sections of $L^{\otimes k}$ are, for large values of $k$, canonical up to
isotopy (this weak uniqueness statement holds in the same sense as that
of Theorem 2.1). Therefore, the topology of the branch curve $D=f(R)$ can
be used to define symplectic invariants, provided that one takes into
account the possibility of cancellations or creations of pairs of nodes
with opposite orientations in isotopies of branched coverings.

Most of the results cited below were obtained in
a joint work with L.~Katzarkov \cite{AK}.

\subsection{Quasiholomorphic maps to $\CP^2$}

In order to study the topology of the singular plane curve $D$, it is
natural to try to adapt the braid group techniques previously used by 
Moishezon and Teicher in the algebraic case (see e.g.\ \cite{moi1},
\cite{moi2}, \cite{teicher}). However, in order to apply
this method it is necessary to ensure that the branch curve satisfies
suitable transversality properties with respect to a generic projection
map from $\CP^2$ to $\CP^1$. This leads naturally to the notion of
{\it quasiholomorphic covering} introduced in \cite{AK}, which we now
describe carefully. 

We slightly rephrase the conditions listed in \cite{AK} in such a way 
that they extend naturally to the higher dimensional case~; the same 
definitions will be used again in \S 4. It is important to be aware that 
these concepts only apply to sequences of objects obtained for increasing 
values of the degree $k$~; the general strategy is always to work
simultaneously with a whole family of sections indexed by the parameter
$k$, in order to ultimately ensure the desired properties for large values
of $k$. We start with the following terminology~:

\begin{defn}
A sequence of sections $s_k$ of complex vector bundles $E_k$ over $X$
(endowed with Hermitian metrics and connections) is {\em asymptotically 
holomorphic} if there exist constants $C_j$ independent of $k$ such that
$|\nabla^j s_k|_{g_k}\le C_j$ and $|\nabla^{j-1}\dbar s_k|_{g_k}\le C_j
k^{-1/2}$ for all $j$, all norms being evaluated with respect to the
rescaled metric $g_k=kg$ on $X$. 

The sections $s_k$ are {\em uniformly transverse to $0$} if there
exists a constant $\gamma>0$ such that, at every point $x\in X$ where
$|s_k(x)|\le\gamma$, the covariant derivative $\nabla s_k(x)$ is surjective 
and has a right inverse of norm less than $\gamma^{-1}$ w.r.t. $g_k$
(we then say that $s_k$ is $\gamma$-transverse to $0$).
\end{defn}

It is easy to check that, if the sections $s_k$ are asymptotically holomorphic
and uniformly transverse to $0$, then for large $k$ their zero sets are
smooth approximately holomorphic symplectic submanifolds.

Also observe that, in the case where the rank of the bundle $E_k$ is greater 
than the dimension of $X$, the surjectivity condition imposed by
transversality is never satisfied~; $\gamma$-transversality to $0$ then
means that the norm of the section is greater than $\gamma$ at every point
of $X$.

\begin{defn}
A sequence of projective maps $f_k:X\to\CP^2$ determined by
asymptotically holomorphic sections $s_k=(s_k^0,s_k^1,s_k^2)$ of
$\C^3\otimes L^{\otimes k}$ for $k\gg 0$ is {\em quasiholomorphic} if there 
exist constants $C_j$, $\gamma$, $\delta$ independent of $k$, 
almost-complex structures $\tilde{J}_k$ on $X$, and finite sets
$\mathcal{C}_k,\mathcal{T}_k,\mathcal{I}_k\subset X$ such that the following 
properties hold (using $\tilde{J}_k$ to define the $\dbar$ operator)~: 

$(0)$ $|\nabla^j(\tilde{J}_k-J)|_{g_k}\le C_j k^{-1/2}$ for every $j\ge 0$~;
$\tilde{J}_k=J$ outside of the $2\delta$-neighborhood of 
$\mathcal{C}_k\cup\mathcal{T}_k\cup\mathcal{I}_k$~; 
$\tilde{J}_k$ is integrable in the $\delta$-neighborhood of
$\mathcal{C}_k\cup\mathcal{T}_k\cup\mathcal{I}_k$~; 

$(1)$ the section $s_k$ of $\C^3\otimes L^{\otimes k}$ is
$\gamma$-transverse to $0$~;

$(2)$ $|\nabla f_k(x)|_{g_k}\ge \gamma$ at every point $x\in X$~; 

$(3)$ the $(2,0)$-Jacobian $\mathrm{Jac}(f_k)=\bigwedge^2\partial f_k$ 
is $\gamma$-transverse to $0$~; in particular it vanishes transversely 
along a smooth symplectic curve $R_k\subset X$ (the branch curve).

$(3')$ the restriction of $\dbar f_k$ to
$\mathrm{Ker}\,\partial f_k$ vanishes at every point of $R_k$~;

$(4)$ the quantity $\partial(f_{k|R_k})$, which can be seen as a section
of a line bundle over $R_k$, is $\gamma$-transverse to $0$ and vanishes at 
the finite set $\mathcal{C}_k$ (the cusp points of $f_k$)~; in particular
$f_k(R_k)=D_k$ is an immersed symplectic curve away from the image of 
$\mathcal{C}_k$~;

$(4')$ $f_k$ is $\tilde{J}_k$-holomorphic over the $\delta$-neighborhood of
$\mathcal{C}_k$~;

$(5)$ the section $(s_k^0,s_k^1)$ of $\C^2\otimes L^{\otimes k}$ is
$\gamma$-transverse to $0$~; as a consequence $D_k$ remains away from the
point $(0:0:1)$~;

$(6)$ let $\pi:\CP^2-\{(0\!:\!0\!:\!1)\}\to\CP^1$ be the map defined by
$\pi(x\!:\!y\!:\!z)=(x\!:\!y)$, and let $\phi_k=\pi\circ f_k$. Then the 
quantity $\partial(\phi_{k|R_k})$ is $\gamma$-transverse to $0$ over $R_k$, 
and it vanishes over the union of $\mathcal{C}_k$ with the finite set 
$\mathcal{T}_k$ (the tangency points of the branch curve $D_k$ with respect
to the projection $\pi$)~; 

$(6')$ $f_k$ is $\tilde{J}_k$-holomorphic over the $\delta$-neighborhood of
$\mathcal{T}_k$~;

$(7)$ the projection $f_k:R_k\to D_k$ is injective outside of the singular
points of $D_k$, and the self-intersections of $D_k$ are transverse double
points. Moreover, all special points of $D_k$ (cusps, nodes, tangencies)
lie in different fibers of the projection $\pi$, and none of them lies
in $\pi^{-1}(0\!:\!1)$~;

$(8)$ the section $s_k^0$ of $L^{\otimes k}$ is $\gamma$-transverse to $0$~;

$(8')$ $R_k$ intersects the zero set of $s_k^0$ at the points of 
$\mathcal{I}_k$~; $f_k$ is 
$\tilde{J}_k$-holomorphic over the $\delta$-neighborhood of
$\mathcal{I}_k$.
\end{defn}

\begin{rem}
Definition 3.2 is slightly stronger than the definition
given in \cite{AK}. Most notably, property $(8)$, which ensures that the 
fiber of $\pi\circ f_k$ above $(0\!:\!1)$ enjoys suitable genericity 
properties, has been added for our purposes. Similarly, condition $(6')$ 
is significantly stronger than in \cite{AK}, where it was only required 
that $\dbar f_k$ vanish at the points of $\mathcal{T}_k$. These extra 
conditions only require minor modifications of the arguments, while 
allowing the inductive construction described in \S 5 to be largely
simplified.
\end{rem}

Observe that, because of property $(0)$, the notions of asymptotic
holomorphicity with respect to $J$ or $\tilde{J}_k$ coincide. Moreover,
even though $\tilde{J}_k$ is used implicitly thoughout the definition, 
the choice of $J$ or $\tilde{J}_k$ is irrelevant as far as transversality
properties are concerned since they differ by $O(k^{-1/2})$.

Property $(1)$ means that $s_k$ is everywhere bounded from below by 
$\gamma$~; this implies that the projective map $f_k$ is well-defined, and
that $|\nabla^j f_k|_{g_k}=O(1)$ and $|\nabla^{j-1}\dbar f_k|_{g_k}=
O(k^{-1/2})$ for all $j$. The second property can be interpreted in terms of
transversality to the codimension 4 submanifold in the space of $1$-jets
given by the equation $\partial f=0$. Properties $(3)$ and $(3')$ yield the
correct structure near generic points of the branch curve~: the transverse
vanishing of $\mathrm{Jac}(f_k)$ implies that the branching order is $2$, 
and the compatibility property $(3')$ ensures that $\dbar f_k$ remains much 
smaller than $\partial f_k$ in all directions, which is needed to obtain the 
correct local model.

Properties $(4)$ and $(4')$ determine the structure of the covering near the
cusp points. More precisely, observe that along $R_k$ the tangent plane
field $TR_k$ and the plane field $\mathrm{Ker}\,\partial f_k$ coincide
exactly at the cusp points~; condition $(4)$ expresses that these two
plane fields are transverse to each other (in \cite{A2} and \cite{AK} this
condition was formulated in terms of a more complicated quantity; the two
formulations are easily seen to be equivalent). This implies that cusp
points are isolated and non-degenerate. The compatibility condition $(4')$ 
then ensures that the expected local model indeed holds. 

The remaining conditions are used to ensure the compatibility of the branch
curve $D_k=f_k(R_k)$ with the projection $\pi$ to $\CP^1$. In particular,
the transversality condition $(6)$ and the corresponding compatibility
condition $(6')$ imply that the points where the branch curve $D_k$ fails
to be transverse to the fibers of $\pi$ are isolated non-degenerate tangency
points. Moreover, property $(7)$ states that the curve $D_k$ is transverse
to itself. This implies that $D_k$ is a {\it braided curve} in the
following sense~:

\begin{defn}
A real $2$-dimensional singular submanifold $D\subset\CP^2$ is a {\em
braided curve} if it satisfies the following properties~: $(1)$ the only 
singular points of $D$ are cusps (with positive orientation) and transverse
double points (with either orientation)~; $(2)$ the point $(0:0:1)$ does not
belong to $D$~; $(3)$ the fibers of the projection $\pi:(x:y:z)\mapsto (x:y)$
are everywhere transverse to $D$, except at a finite set of nondegenerate
tangency points where a local model for $D$ in orientation-preserving
coordinates is $z_2^2=z_1$~; $(4)$ the cusps, nodes and tangency points are
all distinct and lie in different fibers of $\pi$.
\end{defn}

We will see in \S 3.2 that these properties are precisely those needed
in order to apply the braid monodromy techniques of Moishezon-Teicher to
the branch curve $D_k$.

The main result of \cite{AK} can be formulated as follows~:

\begin{thm}[\cite{A2},\cite{AK}]
For $k\gg 0$, it is possible to find asymptotically holomorphic sections of 
$\C^3\otimes L^{\otimes k}$ such that the corresponding projective maps 
$f_k:X\to\CP^2$ are quasiholomorphic branched coverings. Moreover, for large
$k$ these coverings are canonical up to isotopy and up to cancellations of
pairs of nodes in the branch curves $D_k$.
\end{thm}

The uniqueness statement is to be understood in the same weak sense as
for Theorem 2.1~: given two sequences of quasiholomorphic branched coverings
(possibly for different choices of almost-complex structures on $X$), for
large $k$ it is possible to find an interpolating one-parameter family of
quasiholomorphic coverings, the only possible non-trivial phenomenon being
the cancellation or creation of pairs of nodes in the branch curve for
certain parameter values.

The proof of Theorem 3.1 follows a standard pattern~: in order to construct
quasiholomorphic coverings, one starts with any sequence of asymptotically
holomorphic sections of $\C^3\otimes L^{\otimes k}$ and proceeds by
successive perturbations in order to obtain all the required properties,
starting with uniform transversality. Since transversality is an open
condition, it is preserved by the subsequent perturbations. 

So the first 
part of the proof consists in obtaining, by successive perturbation 
arguments, the transversality properties $(1)$, $(2)$, $(3)$ and $(4)$
of Definition 3.2 as in \cite{A2}, $(5)$ and $(6)$ as in \cite{AK}, and also
$(8)$ by a direct application of the result of \cite{D1}. The argument is
notably more technical in the case of $(4)$ and $(6)$ because the 
transversality conditions involve derivatives along the branch curve, but
these can actually all be thought as immediate applications of the general
transversality principle mentioned in the Introduction.

The second part of the proof, which is comparatively easier, deals with
the compatibility conditions. The idea is to ensure these properties by
perturbing the sections $s_k$ by quantities bounded by $O(k^{-1/2})$, which
clearly affects neither holomorphicity nor transversality properties. 
One first chooses suitable almost-complex structures $\tilde{J}_k$ differing
from $J$ by $O(k^{-1/2})$ and integrable near the finite set
$\mathcal{C}_k\cup\mathcal{T}_k\cup\mathcal{I}_k$. It is then possible to
perturb $f_k$ near these points in order to obtain
conditions $(4')$, $(6')$ and $(8')$, by the same argument as in \S 4.1
of \cite{A2}. Next, a generic small perturbation 
yields the self-transversality of $D$ (property $(7)$).
Finally, a suitable perturbation yields property $(3')$ along the branch 
curve without modifying $R_k$ and $D_k$ and without affecting the 
other compatibility properties.

The uniqueness statement is obtained by showing that, provided that $k$ is 
large enough, all the arguments extend verbatim to one-parameter families 
of sections. Therefore, given two sequences of quasiholomorphic coverings,
one starts with a one-parameter family of sections interpolating between
them in a trivial way and perturbs it in such a way that the
required properties hold for all parameter values (with the exception of
$(7)$ when a node cancellation occurs). Since this construction can be
performed in such a way that the two end points of the one-parameter
family are not affected by the perturbation, the isotopy result follows
immediately.

The reader is referred to \cite{A2} and \cite{AK} for more details
(incorporating requirement $(8)$ in the arguments is a trivial task).
 
\subsection{Braid monodromy invariants}
We now describe the monodromy invariants that naturally arise from the
quasiholomorphic coverings described in the previous section. This is a
relatively direct extension to the symplectic framework of the braid
group techniques studied by Moishezon and Teicher in the algebraic case
(see \cite{moi1}, \cite{moi2}, \cite{teicher}).

Recall that the braid group on $d$ strings is the fundamental group 
$B_d=\pi_1(\mathcal{X}_d)$ of the space $\mathcal{X}_d$ of unordered 
configurations of $d$ distinct points in the plane $\R^2$. A braid can
therefore be thought as a motion of $d$ points in the plane. An alternate
description involves compactly supported orientation-preserving
diffeomorphisms of $\R^2$ which globally preserve a set of $d$ given
points~: $B_d=\pi_0(\mathrm{Diff}^+_c(\R^2,\{q_1,\dots,q_d\}))$.
The group $B_d$ is generated by {\it half-twists}, i.e.\ braids in which
two of the $d$ points rotate around each other by 180 degrees while the
other points are preserved. For more details see \cite{birman}.

Consider a braided curve $D\subset\CP^2$ (see Definition 3.3) of fixed 
degree $d$, for example the branch curve of a quasiholomorphic covering as 
given by Theorem 3.1. Projecting to $\CP^1$ via the map $\pi$ makes $D$ a 
singular branched covering of $\CP^1$. The picture is the following~:

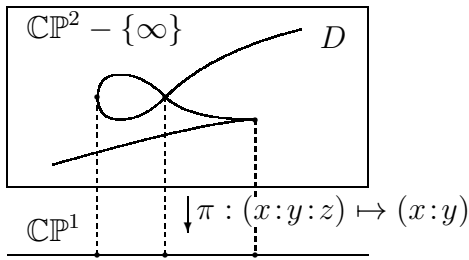
\begin{center}
\setlength{\unitlength}{0.6mm}
\begin{picture}(80,63)(-40,-20)
\put(0,-2){\vector(0,-1){8}}
\put(2,-7){$\pi:(x\!:\!y\!:\!z)\mapsto (x\!:\!y)$}
\put(-40,-15){\line(1,0){80}}
\put(-38,-12){ $\CP^1$}
\put(-40,0){\line(1,0){80}}
\put(-40,40){\line(1,0){80}}
\put(-40,0){\line(0,1){40}}
\put(40,0){\line(0,1){40}}
\put(-38,33){ $\CP^2-\{\infty\}$}
\put(27,31){ $D$}

\multiput(-20,20)(0,-2){18}{\line(0,-1){1}}
\multiput(-5,20)(0,-2){18}{\line(0,-1){1}}
\multiput(15,15)(0,-2){9}{\line(0,-1){1}}
\multiput(15,-9)(0,-2){3}{\line(0,-1){1}}
\put(-20,-15){\circle*{1}}
\put(-5,-15){\circle*{1}}
\put(15,-15){\circle*{1}}

\qbezier(25,35)(5,30)(-5,20)
\qbezier(-5,20)(-10,15)(-15,15)
\qbezier(-15,15)(-20,15)(-20,20)
\qbezier(-20,20)(-20,25)(-15,25)
\qbezier(-15,25)(-10,25)(-5,20)
\qbezier(-5,20)(0,15)(15,15)
\qbezier(15,15)(5,15)(-30,5)
\put(-20,20){\circle*{1}}
\put(-5,20){\circle*{1}}
\put(15,15){\circle*{1}}
\end{picture}
\end{center}

Let $p_1,\dots,p_r$ be the images by $\pi$ of the special points of $D$
(nodes, cusps and tangencies). Observing that the fibers of $\pi$
are complex lines (or equivalently real planes) which generically intersect
$D$ in $d$ points, we easily get that the monodromy of the map $\pi_{|D}$
around the fibers above $p_1,\dots,p_r$ takes values in the braid group
$B_d$.

The monodromy around one of the points $p_1,\dots,p_r$ is as follows.
In the case of a tangency point, a local model for the curve $D$ is
$y^2=x$ (with projection to the $x$ factor), so one easily checks that
the monodromy is a half-twist exchanging two sheets of $\pi_{|D}$. 
Since all half-twists in $B_d$ are conjugate, it is possible to write
this monodromy in the form $Q^{-1}X_1Q$, where $Q\in B_d$ is any braid
and $X_1$ is a fixed half-twist (aligning the points $q_1,\dots,q_d$ in 
that order along the real axis, $X_1$ is the half-twist exchanging the 
points $q_1$ and $q_2$ along a straight line segment). In the case of a
transverse double point with positive intersection, the local model
$y^2=x^2$ implies that the monodromy is the square of a half-twist, which
can be written in the form $Q^{-1}X_1^2Q$. The monodromy around
a double point with negative intersection is the mirror image of the
previous case, and can therefore be written as $Q^{-1}X_1^{-2}Q$.
Finally, the monodromy around a cusp (local model $y^2=x^3$) is the cube
of a half-twist and can be expressed as $Q^{-1}X_1^3Q$.

However, in order to describe the monodromy 
automorphisms as braids, one needs to identify up to compactly supported
diffeomorphisms the fibers of $\pi$ with a reference plane $\R^2$. This 
implicitly requires a trivialization of the 
fibration $\pi$, which is not available over all of $\CP^1$. Therefore,
as in the case of Lefschetz pencils, it is necessary to restrict oneself
to the preimage of an affine subset $\C\subset\CP^1$, by removing 
the fiber above the point at infinity (which may easily be assumed to be
regular). So the monodromy map is only defined as a group homomorphism
\begin{equation}
\rho:\pi_1(\C-\{p_1,\dots,p_r\})\to B_d.
\end{equation}
Since the fibration $\pi$ defines a line bundle of degree $1$ over $\CP^1$,
the monodromy around the fiber at infinity is given by the {\it full twist}
$\Delta^2$, i.e.\ the braid which corresponds to a rotation of all points 
by 360 degrees ($\Delta^2$ generates the center of $B_d$).

Therefore, choosing as in \S 2 a system of generating loops in
$\C-\{p_1,\dots,p_r\}$, we can express the monodromy by a factorization of
$\Delta^2$ in the braid group~:
\begin{equation}
\Delta^2=\prod_{j=1}^r Q_j^{-1}X_1^{r_j}Q_j,
\end{equation}
where the elements $Q_j\in B_d$ are arbitrary braids and the degrees 
$r_j\in\{1,\pm 2,3\}$ depend on the types of the special points lying 
above $p_j$.

As in the case of Lefschetz pencils, this braid factorization, which
completely characterizes the braided curve $D$ up to isotopy, is only
well-defined up to two algebraic operations: simultaneous conjugation of
all factors by a given braid in $B_d$, and Hurwitz moves. As previously,
simultaneous conjugation reflects the different possible choices of 
an identification diffeomorphism between the fiber of $\pi$ above the base 
point and the standard plane $(\C,\{q_1,\dots,q_d\})$, while Hurwitz moves
arise from changes in the choice of a generating system of loops in
$\C-\{p_1,\dots,p_r\}$.

Starting with any braid factorization of the form (4), it is possible to 
reconstruct a braided curve $D$ in a canonical way up to isotopy (see
\cite{AK}; similar statements were also obtained by Moishezon, Teicher
and Catanese).
Moreover, one easily checks that factorizations which differ only by global 
conjugations and Hurwitz moves lead to isotopic braided curves (each such 
operation amounts to a diffeomorphism isotopic to the identity, obtained in 
the case of a Hurwitz move by lifting by $\pi$ a diffeomorphism of $\CP^1$, 
and in the case of a global conjugation by a diffeomorphism in each of the 
fibers of $\pi$).

Moreover, it is important to observe that every braided curve $D$ can be
made symplectic by a suitable isotopy. In fact, it is sufficient to perform
a radial contraction in all the fibers of $\pi$, which brings the given
curve into an arbitrarily small neighborhood of the zero section of $\pi$
(the complex line $\{z=0\}$ in $\CP^2$). The tangent space to $D$ is then
very close to that of the complex line (and therefore symplectic) everywhere
except near the tangency points; verifying that the property also holds near
tangencies by means of the local model, one obtains that $D$ is symplectic.
\medskip

We now briefly describe the structure of the fundamental group
$\pi_1(\CP^2-D)$. Consider a generic fiber of $\pi$, intersecting $D$ in
$d$ points $q_1,\dots,q_d$. Then the inclusion map $i:\C-\{q_1,\dots,q_d\}
\to\CP^2-D$ induces a surjective homomorphism on fundamental groups.
Therefore, a generating system of loops $\gamma_1,\dots,\gamma_d$ in
$\C-\{q_1,\dots,q_d\}$ provides a set of generators for $\pi_1(\CP^2-D)$
({\it geometric generators}). Because the fiber of $\pi$ can be compactified
by adding the pole of the projection, an obvious relation is 
$\gamma_1\dots\gamma_d=1$. Moreover, each special point of the curve $D$,
or equivalently every term in the braid factorization, determines a relation 
in $\pi_1(\CP^2-D)$ in a very explicit way.

Namely, recall that there exists a natural right action of $B_d$ on the 
free group $F_d=\pi_1(\C-\{q_1,\dots,q_d\})$, that we shall denote by $*$, 
and consider a factor $Q_j^{-1}X_1^{r_j}Q_j$ in (4). Then, if $r_j=1$, the 
tangency point above $p_j$ yields the relation 
$\gamma_1 * Q_j=\gamma_2 * Q_j$ (the two elements 
$\gamma_1 * Q_j$ and $\gamma_2 * Q_j$ correspond to small loops going 
around the two sheets of $\pi_{|D}$ that merge at the tangency 
point). Similarly, in the case of a node ($r_j=\pm 2$), the relation is
$[\gamma_1 * Q_j,\gamma_2 * Q_j]=1$. Finally, in the case of a cusp
($r_j=3$), the relation becomes $(\gamma_1\gamma_2\gamma_1)*Q_j=
(\gamma_2\gamma_1\gamma_2)*Q_j$.
It is a classical result that $\pi_1(\CP^2-D)$ is exactly the quotient of
$F_d=\langle\gamma_1,\dots,\gamma_d\rangle$ by the above-listed relations.
\medskip

Given a branched covering map $f:X\to\CP^2$ with branch curve $D$, it is
easy to see that the topology of $X$ is determined by a group homomorphism
from $\pi_1(\CP^2-D)$ to the symmetric group $S_n$ of order $n=\deg f$.
Considering a generic fiber of $\pi$ which intersects $D$ in $d$ points
$q_1,\dots,q_d$, the restriction of $f$ to
its preimage $\Sigma$ is a $n$-sheeted branched covering map from 
$\Sigma$ to $\C$ with branch points $q_1,\dots,q_d$. This covering is 
naturally described by a monodromy representation
\begin{equation}
\theta:\pi_1(\C-\{q_1,\dots,q_d\})\to S_n.
\end{equation}
Because the branching index is $2$ at a generic point of the branch curve 
of $f$, the group homomorphism $\theta$ maps geometric generators to 
transpositions. Also, $\theta$ necessarily factors through the surjective 
homomorphism $i_*:\pi_1(\C-\{q_1,\dots,q_d\})\to\pi_1(\CP^2-D)$, because 
the covering $f$ is defined everywhere, and the resulting map from
$\pi_1(\CP^2-D)$ to $S_n$ is exactly what is needed to recover the
4-manifold $X$ from the branch curve $D$. The properties of $\theta$
are summarized in the following definition due to Moishezon~:

\begin{defn} A geometric monodromy representation associated to a
braided curve $D\subset\CP^2$ is a surjective group homomorphism
$\theta$ from the free group $\pi_1(\C-\{q_1,\dots,q_d\})=F_d$ to the 
symmetric group $S_n$ of order $n$, mapping the geometric generators
$\gamma_i$ (and thus also the $\gamma_i*Q_j$) to transpositions, and such
that\medskip

$\theta(\gamma_1\dots\gamma_d)=1,$

$\theta(\gamma_1*Q_j)=\theta(\gamma_2*Q_j)$ if $r_j=1$,

$\theta(\gamma_1*Q_j)$ and $\theta(\gamma_2*Q_j)$ are distinct and
commute if $r_j=\pm 2$,

$\theta(\gamma_1*Q_j)$ and $\theta(\gamma_2*Q_j)$ do not commute 
if $r_j=3$.
\end{defn}

Observe that, when the braid factorization defining $D$ is affected
by a Hurwitz move, $\theta$ remains unchanged and the compatibility 
conditions are preserved. On the contrary, when the braid 
factorization is modified by simultaneously conjugating all factors 
by a certain braid $Q\in B_d$, the system of geometric generators
$\gamma_1,\dots,\gamma_d$ changes accordingly, and so the geometric monodromy 
representation $\theta$ should be replaced by $\theta\circ Q_*$, where
$Q_*$ is the automorphism of $F_d$ induced by the braid $Q$.

One easily checks that, given a braided curve $D\subset\CP^2$ and a 
compatible monodromy representation $\theta:F_d\to S_n$, it is possible
to recover a compact 4-manifold $X$ and a branched covering map 
$f:X\to\CP^2$ in a canonical way. Moreover, as observed above we can
assume that the curve $D$ is symplectic; in that case, the branched covering
map makes it possible to endow $X$ with a symplectic structure, canonically
up to symplectic isotopy (see \cite{A2},\cite{AK}~; a similar result has 
also been obtained by Catanese).
\medskip

The above discussion leads naturally to the definition of symplectic
invariants arising from the quasiholomorphic coverings constructed in
Theorem 3.1. However, things are complicated by the fact that the branch 
curves of these coverings are only canonical up to cancellations of double 
points.

On the level of the braid factorization, a pair cancellation amounts to
removing two consecutive factors which are the inverse of each other
(necessarily one must have degree $2$ and the other degree $-2$); the
geometric monodromy representation is not affected. The opposite operation
is the creation of a pair of nodes, in which two factors
$(Q^{-1}X_1^{-2}Q).(Q^{-1}X_1^2 Q)$ are added anywhere in the
factorization~; it is allowed only if the new factorization remains
compatible with the monodromy representation $\theta$, i.e.\ if
$\theta(\gamma_1*Q)$ and $\theta(\gamma_2*Q)$ are commuting disjoint
transpositions.

\begin{defn}
Two braid factorizations (along with the corresponding geometric monodromy 
representations) are m-equivalent if there exists a sequence of operations 
which turns one into the other, each operation being either a global 
conjugation, a Hurwitz move, or a pair cancellation or creation. 
\end{defn}

In conclusion, we get the following result~:

\begin{thm}[\cite{AK}]
The braid factorizations and geometric monodromy representations 
associated to the quasiholomorphic coverings obtained in Theorem 3.1
are, for $k\gg 0$, canonical up to m-equivalence, and define symplectic
invariants of $(X^4,\omega)$.

Conversely, the data consisting of a braid factorization and a geometric
monodromy representation, or a m-equivalence class of such data, 
determines a symplectic 4-manifold in a canonical way up to 
symplectomorphism.
\end{thm}

\subsection{The braid group and the mapping class group}

Let $f:X\to\CP^2$ be a branched covering map, and let $D\subset\CP^2$ be
its branch curve. It is a simple observation that, if $D$ is braided, then
the map $\pi\circ f$ with values in $\CP^1$ obtained by forgetting one of
the components of $f$ topologically defines a Lefschetz pencil. This
pencil is obtained by lifting via the covering $f$ the pencil of lines 
on $\CP^2$ defined by $\pi$, and its base points are the preimages by $f$ 
of the pole of the projection $\pi$.

Moreover, if one starts with the quasiholomorphic coverings given by
Theorem 3.1, then the corresponding Lefschetz pencils coincide for $k\gg 0$
with those obtained by Donaldson in \cite{D2} and described in \S 2.

As a consequence, in the case of a 4-manifold, the invariants described in 
\S 3.2 (braid factorization and geometric monodromy representation) 
completely determine those described in \S 2 (factorizations in mapping 
class groups). It is therefore natural to look for a more explicit
description of the relation between branched coverings and
Lefschetz pencils. This description involves the group of
{\it liftable braids}, which has been studied in a special case by
Birman and Wajnryb in \cite{bw}. We recall the following construction
from \S 5 of \cite{AK}.
\medskip

Let $\mathcal{C}_n(q_1,\dots,q_d)$ be the (finite) set of
all surjective group homomorphisms $F_d\to S_n$ which map each of the
geometric generators $\gamma_1,\dots,\gamma_d$ of $F_d$ to a transposition 
and map their product $\gamma_1\cdots\gamma_d$ to the identity element in 
$S_n$. Each element of $\mathcal{C}_n(q_1,\dots,q_d)$ determines a
simple $n$-fold covering of $\CP^1$ branched at $q_1,\dots,q_d$.

Let $\mathcal{X}_d$ be the space of configurations of $d$ distinct points in
the plane. The set of all simple $n$-fold coverings of $\CP^1$ with $d$
branch points and such that no branching occurs above the point at infinity
can be thought of as a covering $\tilde{\mathcal{X}}_{d,n}$ above 
$\mathcal{X}_d$, in which the fiber above the configuration 
$\{q_1,\dots,q_d\}$ identifies with
$\mathcal{C}_n(q_1,\dots,q_d)$. Therefore, the 
braid group $B_d=\pi_1(\mathcal{X}_d)$ acts on the fiber
$\mathcal{C}_n(q_1,\dots,q_d)$ by deck transformations of the covering
$\tilde{\mathcal{X}}_{d,n}$. In fact, the action of a braid $Q\in B_d$ on
$\mathcal{C}_n(q_1,\dots,q_d)$ is given by $\theta\mapsto \theta\circ Q_*$,
where $Q_*\in\mathrm{Aut}(F_d)$ is the automorphism induced by $Q$ on the
fundamental group of $\C-\{q_1,\dots,q_d\}$.

Fix a base point $\{q_1,\dots,q_d\}$ in $\mathcal{X}_d$, and consider an
element $\theta$ of $\mathcal{C}_n(q_1,\dots,q_d)$ (i.e., a monodromy 
representation $\theta:F_d\to S_n$). Let $p_\theta$ be the corresponding
point in $\tilde{\mathcal{X}}_{d,n}$.

\begin{defn}
The subgroup $B_d^0(\theta)$ of liftable braids is the set of
all the loops in $\mathcal{X}_d$ whose lift at the point $p_\theta$
is a closed loop in $\tilde{\mathcal{X}}_{d,n}$. Equivalently, $B_d^0(\theta)$
is the set of all braids which act on $F_d=\pi_1(\C-\{q_1,\dots,q_d\})$ in 
a manner compatible with the covering structure defined by $\theta$.
\end{defn}

In other words, $B_d^0(\theta)$ is the set of all braids $Q$ such that
$\theta\circ Q_*=\theta$, i.e.\ the stabilizer of $\theta$ with respect
to the action of $B_d$ on $\mathcal{C}_n(q_1,\dots,q_d)$.

There exists a natural bundle $\mathcal{Y}_{d,n}$ over 
$\tilde{\mathcal{X}}_{d,n}$ (the {\it universal curve}) whose fiber is a 
Riemann surface of genus $g=1-n+(d/2)$ with $n$ marked points.
Each of these Riemann surfaces naturally carries a structure of branched 
covering of $\CP^1$, and the marked points are the preimages of the
point at infinity.

Given an element $Q$ of $B_d^0(\theta)\subset B_d$, it can be lifted to
$\tilde{\mathcal{X}}_{d,n}$ as a loop based at the point $p_\theta$, and
the monodromy of the fibration $\mathcal{Y}_{d,n}$ along this loop defines
an element of $\mathrm{Map}_{g,n}$ (the mapping class group of a Riemann 
surface of genus $g$ with $n$ boundary components), which 
we call $\theta_*(Q)$. This defines a group homomorphism
$\theta_*:B_d^0(\theta)\to\mathrm{Map}_{g,n}$.

More geometrically, viewing $Q$ as a compactly supported diffeomorphism of 
the plane preserving $\{q_1,\dots,q_d\}$, the fact that $Q$ belongs to
$B_d^0(\theta)$ means that it can be lifted via the covering map
$\Sigma_g\to\CP^1$ to a diffeomorphism of $\Sigma_g$~; the corresponding
element in the mapping class group is $\theta_*(Q)$.
\medskip

It is easy to check that, when the given monodromy representation $\theta$
is compatible with a braided curve $D\subset\CP^2$, the image of the braid 
monodromy homomorphism $\rho:\pi_1(\C-\{p_1,\dots,p_r\})\to B_d$ describing
$D$ is entirely contained in $B^0_d(\theta)$~: this is because the geometric 
monodromy representation $\theta$ factors through $\pi_1(\CP^2-D)$, on which
the braids in $\mathrm{Im}\,\rho$ act trivially.
Therefore, we can take the image of the braid factorization describing $D$ 
by $\theta_*$ and obtain a factorization in the mapping class group 
$\mathrm{Map}_{g,n}$. One easily checks that $\theta_*(\Delta^2)$ is, as
expected, the twist $\delta_Z$ around the $n$ marked points.

As observed in \cite{AK}, all the factors of degree $\pm 2$ or $3$ in the
braid factorization lie in the kernel of $\theta_*$~; therefore, the only terms
whose contribution to the mapping class group factorization is non-trivial 
are those arising from the tangency points of the branch curve $D$, and each
of these is a Dehn twist. More precisely, the image in $\mathrm{Map}_{g,n}$ 
of a half-twist $Q\in B^0_d(\theta)$ can be constructed as follows. 
Call $\gamma$ the path joining two of the branch points naturally 
associated to the half-twist $Q$
(i.e.\, the path along which the twisting occurs). Among the $n$ lifts of 
$\gamma$ to $\Sigma_g$, only two hit the branch points of the covering~; 
these two lifts have common end points, and together they define a loop 
$\delta$ in $\Sigma_g$. Then the element $\theta_*(Q)$ in 
$\mathrm{Map}_{g,n}$ is the positive Dehn twist along
the loop $\delta$ (see Proposition 4 of \cite{AK}).

In conclusion, the following result holds~:

\begin{prop}
Let $f:X\to\CP^2$ be a branched covering, and assume that its branch curve 
$D$ is braided. Let $\rho:\pi_1(\C-\{p_1,\dots,p_r\})\to B_d^0(\theta)$ and
$\theta:F_d\to S_n$ be the corresponding braid monodromy and geometric
monodromy representation. Then the monodromy map 
$\psi:\pi_1(\C-\{p_1,\dots,p_r\})\to \mathrm{Map}_{g,n}$  of the Lefschetz 
pencil $\pi\circ f$ is given by the identity $\psi=\theta_*\circ\rho$.

In particular, for $k\gg 0$ the symplectic invariants obtained from 
Theorem 2.1 are obtained in this manner from those given by Theorem 3.2.
\end{prop}

\begin{rem}
It is a basic fact that for $n\ge 3$ the group homomorphism
$\theta_*:B_d^0(\theta)\to\mathrm{Map}_{g,n}$ is surjective, and that
for $n\ge 4$ every Dehn twist is the image by $\theta_*$ of a half-twist.
This makes it natural to ask whether every factorization of $\delta_Z$ in
$\mathrm{Map}_{g,n}$ as a product of Dehn twists is the image by $\theta_*$
of a factorization of $\Delta^2$ in $B_d^0(\theta)$ compatible with $\theta$.
This can be reformulated in more geometric terms as the classical problem
of determining whether every Lefschetz pencil is topologically a covering of 
$\CP^2$ branched along a curve with node and cusp singularities
(a similar question replacing pencils by Lefschetz fibrations and $\CP^2$ 
by ruled surfaces also holds~; presently the answer is only known in the
hyperelliptic case, thanks to the results of Fuller, Siebert and Tian).

A natural approach to these problems is to understand the kernel of
$\theta_*$. For example, if one can show that this kernel is generated by
squares and cubes of half-twists (factors of degree $2$ and $3$ compatible
with $\theta$), then the solution naturally follows~: given a decomposition
of $\delta_Z$ as a product of Dehn twists in $\mathrm{Map}_{g,n}$, any lift 
of this word to $B_d^0(\theta)$ as a product of half-twists differs from 
$\Delta^2$ by a product of factors of degree $2$ and $3$ and their inverses.
Adding these factors as needed, one obtains a decomposition of $\Delta^2$ 
into factors of degrees $1$, $\pm 2$ and $\pm 3$~; the branch curve
constructed in this way may have nodes and cusps with reversed orientation,
but it can still be made symplectic. 

Even if the kernel of $\theta_*$ is not generated by factors of degree $2$
and $3$, it remains likely that the result still holds and can be obtained
by starting from a suitable lift to $B_d^0(\theta)$ of the word in 
$\mathrm{Map}_{g,n}$. A better understanding of the structure of 
$\mathrm{Ker}\,\theta_*$ would be extremely useful for this purpose.
\end{rem}

\section{The higher dimensional case}
In this section we extend the results of \S 3 to the case of higher
dimensional symplectic manifolds. In \S 4.1 we prove the existence of
quasiholomorphic maps $X\to\CP^2$ given by triples of sections of
$L^{\otimes k}$ for $k\gg 0$. The topological invariants arising from 
these maps are studied in \S 4.2 and \S 4.3, and the relation with 
Lefschetz pencils is described in \S 4.4.

\subsection{Quasiholomorphic maps to $\CP^2$}

Let $(X^{2n},\omega)$ be a compact symplectic manifold, endowed with a
compatible almost-complex structure $J$. Let $L$ be the same line bundle
as previously (if $\frac{1}{2\pi}[\omega]$ is not integral one works with
a perturbed symplectic form as explained in the introduction). Consider
three approximately holomorphic sections of $L^{\otimes k}$, or equivalently
a section of $\C^3\otimes L^{\otimes k}$. Then the following result states
that exactly the same transversality and compatibility properties can be
expected as in the four-dimensional case~:

\begin{thm}
For $k\gg 0$, it is possible to find asymptotically holomorphic sections of 
$\C^3\otimes L^{\otimes k}$ such that the corresponding $\CP^2$ valued 
projective maps $f_k$ are quasiholomorphic (cf.\ Definition 3.2). 
Moreover, for large $k$ these projective maps are canonical up to 
isotopy and up to cancellations of pairs of nodes in the critical 
curves $D_k$.
\end{thm}

Before sketching a proof of Theorem 4.1, we briefly describe the behavior
of quasiholomorphic maps, which will clarify some of the requirements of
Definition 3.2.

Condition $(1)$ in Definition 3.2 implies that the set $Z_k$ of points where 
the three sections $s_k^0,s_k^1,s_k^2$ vanish simultaneously is a smooth 
codimension 6 symplectic (approximately holomorphic) submanifold. The
projective map $f_k=(s_k^0\!:\!s_k^1\!:\!s_k^2)$ with values in $\CP^2$ is
only defined over the complement of $Z_k$. The behavior near the set of base 
points is similar to what happens for Lefschetz pencils~: in suitable local
approximately holomorphic coordinates, $Z_k$ is given by the equation
$z_1=z_2=z_3=0$, and $f_k$ behaves like the model map
$(z_1,\dots,z_n)\mapsto (z_1\!:\!z_2\!:\!z_3)$. In fact, a map defined
everywhere can be obtained by blowing up $X$ along the submanifold $Z_k$.
The behavior near $Z_k$ being completely specified by condition $(1)$, it
is implicit that all the other conditions on $f_k$ are only to be imposed
outside of a small neighborhood of $Z_k$.

The correct statement of condition $(3)$ of Definition 3.2 in the case
of a manifold of dimension greater than $4$ is a bit tricky. Indeed, 
$\mathrm{Jac}(f_k)=\bigwedge^2\partial f_k$ is a priori a section of the
vector bundle $\Lambda^{2,0}T^*X\otimes f_k^*(\Lambda^{2,0}T\CP^2)$ of
rank $n(n-1)/2$. However, transversality to $0$ in this sense is 
impossible to obtain, as the expected complex codimension of $R_k$ is $n-1$
instead of $n(n-1)/2$. Indeed, the section $\mathrm{Jac}(f_k)$ takes values 
in the non-linear subbundle $\mathrm{Im}(\bigwedge^2)$, whose fibers are of 
dimension $n-1$ at their smooth points (away from the origin). However,
transversality to $0$ does not have any natural definition in this
subbundle, because it is singular along the zero section. 
The problem is very similar to what happens in the construction of 
determinantal submanifolds performed in \cite{presas}. 

In our case, a precise meaning can be given to condition $(3)$ by the
following observation. Near any point $x\in X$, property $(2)$ implies
that it is possible to find
local approximately holomorphic coordinates on $X$ and local complex 
coordinates on $\CP^2$ in which the differential at $x$ of the first 
component of $f_k$ is can be written $\partial f_k^1(x)=\lambda\,dz_1$, 
with $|\lambda|>\gamma/2$. This implies that, near $x$, the projection of
$\bigwedge^2\partial f_k$ to its components along $dz_1\wedge dz_2,\dots,
dz_1\wedge dz_n$ is a quasi-isometric isomorphism. In other words,
the transversality to $0$ of $\mathrm{Jac}(f_k)$ is to be understood as the
transversality to $0$ of its orthogonal projection to the linear subbundle of 
rank $n-1$ generated by $dz_1\wedge dz_2,\dots,dz_1\wedge dz_n$.

Another equivalent approach is to consider the (non-linear) bundle 
$\mathcal{J}^1(X,\CP^2)$ of holomorphic $1$-jets of maps from $X$ to $\CP^2$. 
Inside this bundle, the $1$-jets whose differential is not surjective define 
a subbundle $\Sigma$ of codimension $n-1$, smooth away from the 
stratum $\{\partial f=0\}$. Since this last stratum is avoided by the
$1$-jet of $f_k$ (because of condition $(2)$), the transversality to $0$ of
$\mathrm{Jac}(f_k)$ can be naturally rephrased in terms of estimated
transversality to $\Sigma$ in the bundle of jets (this approach will
be developed in \cite{Atrans}).

With this understood, conditions $(3)$ and $(3')$ imply, as in the 
four-dimensional case, that the set $R_k$ of points where the differential 
of $f_k$ fails to be surjective is a smooth symplectic curve $R_k\subset X$, 
disjoint from $Z_k$, and that the differential of $f_k$ has rank $2$ at every
point of $R_k$. Also, as before, conditions $(4)$ and $(4')$ imply that
$f_k(R_k)=D_k$ is a symplectic curve in $\CP^2$, immersed outside of the
cusp points.

We now describe the proof of Theorem 4.1~; most of the argument is identical
to the 4-dimensional case, and the reader is referred to \cite{A2} and 
\cite{AK} for notations and details.

\begin{proof}[Proof of Theorem 4.1]
The strategy of proof is the same as in the 4-dimensional case. One starts
with an arbitrary sequence of asymptotically holomorphic sections of
$\C^3\otimes L^{\otimes k}$ over $X$, and perturbs it first to obtain the
transversality properties. Provided that $k$ is large enough, each 
transversality property can be obtained over a ball by a small localized 
perturbation, using the local transversality
result of Donaldson (Theorem 12 in \cite{D2}). A globalization argument
then makes it possible to combine these local perturbations into a global
perturbation that ensures transversality everywhere 
(Proposition 3 of \cite{A2}). Since transversality properties are open,
successive perturbations can be used to obtain all the required properties~:
once a transversality property is obtained, subsequent perturbations only
affect it by at most decreasing the transversality estimate.

{\bf Step 1.}
One first obtains the transversality statements in parts $(1)$, $(5)$ and
$(8)$ of Definition 3.2~; as in the 4-dimensional case, these properties
are obtained e.g.\ simply by applying the main result of \cite{A1}.
Observe that all required properties now hold near the base locus $Z_k$ of
$s_k$, so we can assume in the rest of the argument that the points
of $X$ being considered lie away from $Z_k$, and therefore that $f_k$ is 
locally well-defined.

One next ensures condition $(2)$, for which the argument is an immediate
adaptation of that in \S 2.2 of \cite{A2}, the only difference being 
the larger number of coordinate functions.

{\bf Step 2.}
The next property we want to get is condition $(3)$. Here a significant
generalization of the argument in \S 3.1 of \cite{A2} is needed.
The problem reduces, as usual, to showing that the uniform transversality 
to $0$ of $\mathrm{Jac}(f_k)$ can be ensured over a small ball centered at
a given point $x\in X$ by a suitable localized perturbation. As in \cite{A2} 
one can assume that $s_k(x)$ is of the form $(s_k^0(x),0,0)$ and therefore
locally trivialize $\CP^2$ via the quasi-isometric map
$(x\!:\!y\!:\!z)\mapsto (y/x,z/x)$~; this reduces the problem to the study
of a $\C^2$-valued map $h_k$. Because $|\partial f_k|$ is bounded from
below, we can assume (after a suitable rotation) that $|\partial h_k^1(x)|$ 
is greater than some
fixed constant. Also, fixing suitable approximately holomorphic Darboux
coordinates $z_k^1,\dots,z_k^n$ (using Lemma 3 of \cite{A2}, which 
trivially extends to dimensions larger than 4), we can after a rotation
assume that $\partial h_k^1(x)$ is of the form $\lambda\,dz_k^1$, where
the complex number $\lambda$ is bounded from below.

By Lemma 2 of \cite{A2}, there exist asymptotically holomorphic sections
$s_{k,x}^\mathrm{ref}$ of $L^{\otimes k}$ with exponential decay away from
$x$. Define the asymptotically holomorphic 2-forms
$\mu_k^j=\partial h_k^1\wedge \partial(z_k^j s_{k,x}^\mathrm{ref}/s_k^0)$ 
for $2\le j\le n$. At $x$, the 2-form $\mu_k^j$ is proportional to 
$dz_k^1\wedge dz_k^j$~; therefore, over a small neighborhood of $x$, the
transversality to $0$ of $\mathrm{Jac}(f_k)$ in the sense explained above
is equivalent to the transversality to $0$ of the projection of
$\mathrm{Jac}(h_k)$ onto the subspace generated by $\mu_k^2,\dots,\mu_k^n$.
In terms of $1$-jets, the 2-forms $\mu_k^j$ define a local frame in the
normal bundle to the stratum of non-regular maps at $\mathcal{J}^1(f_k)$.
Now, express $\mathrm{Jac}(h_k)$ in the form
$u_k^2\mu_k^2+\dots+u_k^n\mu_k^n+\alpha_k$ over a neighborhood of $x$, where
$u_k^2,\dots,u_k^n$ are complex-valued functions and $\alpha_k$ has no
component along $dz_k^1$. Then, the transversality to $0$ of
$\mathrm{Jac}(f_k)$ is equivalent to that of the $\C^{n-1}$-valued function
$u_k=(u_k^2,\dots,u_k^n)$.

Since the functions $u_k$ are asymptotically holomorphic, using suitable
Darboux coordinates at $x$ we can use Theorem 12 of \cite{D2} to obtain,
for large enough $k$, the existence of constants $w_k^2,\dots,w_k^n$ smaller 
than any given bound $\delta>0$ and
such that $(u_k^2-w_k^2,\dots,u_k^n-w_k^n)$ is $\eta$-transverse to $0$ over
a small ball centered at $x$, where $\eta=\delta(\log\delta^{-1})^{-p}$ ($p$
is a fixed constant). Letting $\tilde{s}_k=(s_k^0,s_k^1,s_k^2-\sum w_k^j
z_k^j s_{k,x}^\mathrm{ref})$ and calling $\tilde{f}_k$ and $\tilde{h}_k$ 
the projective map defined by $\tilde{s}_k$ and the corresponding local 
$\C^2$-valued map, we get that $\mathrm{Jac}(\tilde{h}_k)=\mathrm{Jac}(h_k)-
\sum w_k^j \mu_k^j$, and therefore that $\mathrm{Jac}(\tilde{f}_k)$ is
transverse to $0$ near $x$. Since the perturbation of $s_k$ has exponential
decay away from $x$, we can apply the standard globalization argument to
obtain property $(3)$ everywhere.

{\bf Step 3.}
The next properties that we want to get are $(4)$ and $(6)$. It is possible
to extend the arguments of \cite{A2} and \cite{AK} to the higher dimensional
case~; however this yields a very technical and lengthy argument, so we 
outline here a more efficient strategy following the ideas of \cite{Atrans}. 
Thanks to the previously obtained transversality properties $(1)$ and $(5)$, 
both $f_k$ and $\phi_k$ are well-defined over a neighborhood of $R_k$, so the 
statements of $(4)$ and $(6)$ are well-defined. Moreover, observe that 
property $(6)$ implies property $(4)$, because at any point where 
$\partial(f_{k|R_k})$ vanishes, $\partial(\phi_{k|R_k})$ necessarily
vanishes as well, and if it does so transversely then the same is true for
$\partial(f_{k|R_k})$ as well. So we only focus on $(6)$.

This property can be rephrased in terms of transversality to the codimension
$n$ stratum $S:\{\partial(\phi_{|R})=0\}$ in the bundle 
$\mathcal{J}^2(X,\CP^2)$ of holomorphic 2-jets of maps from $X$ to $\CP^2$.
However this stratum is singular, even away from the substratum $S_{nt}$
corresponding to the non-transverse vanishing of $\mathrm{Jac}(f)$~; in fact
it is reducible and comes as a union $S_1\cup S_2$, where
$S_1:\{\mathrm{Jac}(f)=0,\ \partial(f_{|R})=0\}$ is the stratum corresponding 
to non-immersed points of the branch curve, and $S_2:\{\partial\phi=0\}$ is 
the stratum corresponding to tangency points of the branch curve. Therefore, 
one first needs to ensure transversality with respect to $S_0=S_1\cap
S_2:\{\partial\phi=0,\ \partial(f_{|R})=0\}$, which is a smooth codimension
$n+1$ stratum (``vertical cusp points of the branch curve'') away from
$S_{nt}$.

{\bf Step 3a.}
We first show that a small perturbation can be used to make sure that
the quantity $(\partial\phi_k,\partial(f_{k|R_k}))$ remains bounded from
below, i.e.\ that given any point $x\in X$, either $\partial\phi_k(x)$ is
larger than a fixed constant, or $x$ lies at more than a fixed distance from 
$R_k$, or $x$ lies close to a point of $R_k$ where $\partial(f_{k|R_k}))$ is 
larger than a fixed constant. Since this transversality property is local and
open, we can obtain it by successive small localized perturbations, as for
the previous properties.

Fix a point $x\in X$, and assume that $\partial\phi_k(x)$ is small
(otherwise no perturbation is needed). By property $(5)$, we know that
necessarily $(s_k^0,s_k^1)$ is bounded away from zero at $x$~; a rotation 
in the first two coordinates makes it possible to assume that $s_k^1(x)=0$
and $s_k^0$ is bounded from below near $x$. As above, we replace $f_k$ by
the $\C^2$-valued map $h_k=(h_k^1,h_k^2)$, where $h_k^i=s_k^i/s_k^0$.
By assumption, we get that $\partial h_k^1(x)$ is small. This implies in
particular that $\mathrm{Jac}(f_k)$ is small at $x$, and therefore
property $(3)$ gives a lower bound on its covariant derivative. Moreover,
by property $(2)$ we also have a lower bound on $\partial h_k^2(x)$, which
after a suitable rotation can be assumed equal to $\lambda\,dz_k^1$ for some
$\lambda\neq 0$.
So, as above we can express $\bigwedge^2 \partial f_k$ by looking at its 
components
along $dz_k^1\wedge dz_k^j$ for $2\le j\le n$~; we again define the $2$-forms
$\mu_k^j=\partial h_k^2\wedge \partial(z_k^js_{k,x}^\mathrm{ref}/s_k^0)$,
and the functions $u_2,\dots,u_n$ are defined as previously.
Define a $(n,0)$-form $\theta$ over a neighborhood of $x$ by 
$\theta=\partial u_2\wedge\dots\wedge\partial u_n\wedge\partial
h_k^2$~: at points of $R_k$, the vanishing of $\theta$ is equivalent to
that of $\partial h^2_{k|R_k}$, or equivalently to that of $\partial
f_{k|R_k}$. So our aim is to show that the quantity $(\partial
h_k^1,\theta)$, which is a section of a rank $n+1$ bundle $\mathcal{E}_0$
near $x$, can be made bounded from below by a small perturbation.

For this purpose, we first show the existence of complex-valued polynomials
$(P^1_j,P^2_j)$ and local sections $\epsilon_j$ of $\mathcal{E}_0$, 
$1\le j\le n+1$, such that~: 

(a) for any coefficients $w_j\in\C$, replacing the 
given sections of $L^{\otimes k}$ by $(s_k^0,s_k^1+\sum w_j P^1_j 
s_{k,x}^\mathrm{ref},s_k^2+\sum w_j P^2_j s_{k,x}^\mathrm{ref})$ affects 
$(\partial h_k^1,\theta)$ by the addition of $\sum w_j \epsilon_j +
O(w_j^2)$~; 

(b) the sections $\epsilon_j$ define a local frame in
$\mathcal{E}_0$, and $\epsilon_1\wedge\dots\wedge\epsilon_{n+1}$ is bounded
from below by a universal constant.

First observe that, by property $(3)$, $\partial u_2\wedge\dots\wedge
\partial u_n$ is bounded from below near $x$, whereas we may assume that
$\theta=\partial u_2\wedge\dots\wedge \partial u_n\wedge\partial h_k^2$
is small (otherwise no perturbation is needed). Therefore,
$\partial h_k^2$ (which at $x$ is colinear to $dz_k^1$) lies close to the 
span of the $\partial u_j$.
In particular, after a suitable rotation in the $n-1$ last coordinates on
$X$, we can assume that $\partial u_2\wedge \partial h_k^2$ is small at $x$.
On the other hand, we know that there exists $j_0\neq 1$ such that
$dz_k^{j_0}$ lies far from the span of the $\partial u_j(x)$. We then define
$P^1_{n+1}=z_k^2z_k^{j_0}$ and $P^2_{n+1}=0$. Adding to $s_k^1$ a quantity
of the form $w\,z_k^2z_k^{j_0}\,s_{k,x}^\mathrm{ref}$ does not affect 
$\partial h_k(x)$, but affects $\partial u_2(x)$ by the addition of 
a non-trivial multiple of $dz_k^{j_0}$, and similarly affects $\partial
u_{j_0}(x)$ by the addition of a non-trivial multiple of $dz_k^2$. The other
$\partial u_j(x)$ are not affected. Therefore, $\theta(x)$ changes by an
amount of
$$cw\,dz_k^{j_0}\wedge\partial u_3\wedge\dots\wedge\partial u_n\wedge 
\partial h_k^2+c'w\,\partial u_2\wedge\dots\wedge dz_k^2\wedge\dots\wedge 
\partial u_n\wedge \partial h_k^2+O(w^2),$$
where the constants $c$ and $c'$ are bounded from above and below. The first
term is bounded from below by construction, while the second term is only
present if $j_0\neq 2$ (this requires $n\ge 3$), and in that case it is
small because $\partial u_2\wedge\partial h_k^2$ is small. Therefore, the
local section $\epsilon_{n+1}$ of $\mathcal{E}_0$ naturally corresponding
to such a perturbation is of the form $(0,\epsilon'_{n+1})$ at $x$, where
$\epsilon'_{n+1}$ is bounded from below.

Next, for $1\le j\le n$ we define $P^1_j=z_k^j$ and $P^2_j=0$, and observe
that adding $w\,z_k^j\,s_{k,x}^\mathrm{ref}$ to $s_k^1$ affects 
$\partial h_k^1(x)$ by adding a nontrivial multiple of $dz_k^j$. Therefore,
the local section of $\mathcal{E}_0$ corresponding to this perturbation is
at $x$ of the form $\epsilon_j(x)=(c''dz_k^j,\epsilon'_j)$, where $c''$ 
is a constant bounded from below.

It follows from this argument that the chosen perturbations $P^1_j$ and
$P^2_j$ for $1\le j\le n+1$, and the corresponding local sections
$\epsilon_j$ of $\mathcal{E}_0$, satisfy the conditions (a) and (b) expressed
above. Observe that, because $\epsilon_j$ define a local frame at $x$
and $\epsilon_1\wedge\dots\wedge \epsilon_{n+1}$ is bounded from below at
$x$, the same properties remain true over a ball of fixed radius
around $x$.

Now that a local approximately holomorphic frame in $\mathcal{E}_0$ is
given, we can write $(\partial h_k^1,\theta)$ in the form 
$\sum \zeta_j \epsilon_j$
for some complex-valued functions $\zeta_j$~; it is easy to check that
these functions are asymptotically holomorphic. Therefore, we can again
use Theorem 12 of \cite{D2} to obtain, if $k$ is large enough, the existence 
of constants $w_1,\dots,w_{n+1}$ smaller than any given bound $\delta>0$ 
and such that $(\zeta_1-w_1,\dots,\zeta_{n+1}-w_{n+1})$ is bounded from
below by $\eta=\delta(\log\delta^{-1})^{-p}$ ($p$ is a fixed constant)
over a small ball centered at $x$. Letting $\tilde{s}_k=(s_k^0,
s_k^1-\sum w_j P^1_j s_{k,x}^\mathrm{ref},s_k^2-\sum w_j P^2_j 
s_{k,x}^\mathrm{ref})$ and calling $\tilde{f}_k$, $\tilde{h}_k$ and
$\tilde\theta$ the projective map defined by $\tilde{s}_k$ and the 
corresponding local maps, we get that $(\partial\tilde{h}^1_k,\tilde\theta)$ 
is by construction bounded from below by $c_0\eta$, for a fixed constant 
$c_0$~; indeed, observe that the non-linear term $O(w^2)$ in the perturbation 
formula does not play any significant role, as it is at most of the order 
of $\delta^2\ll\eta$. Since the perturbation of $s_k$ has exponential
decay away from $x$, we can apply the standard globalization argument to
obtain uniform transversality to the stratum
$S_0\subset\mathcal{J}^2(X,\CP^2)$ everywhere.

{\bf Step 3b.}
We now obtain uniform transversality to the stratum $S:\{\mathrm{Jac}(f)=0,
\ \partial(\phi_{|R})=0\}$. The strategy and notations are
the same as above. We again fix a point $x\in X$, and assume that $x$ lies
close to a point of $R_k$ where $\partial(\phi_{k|R_k})$ is small
(otherwise, no perturbation is needed). As above, we can assume that
$s_k^0(x)$ is bounded from below and define a $\C^2$-valued map $h_k$.
Two cases can occur~: either $\partial h_k^1(x)$ is bounded away from zero, 
or it is small and in that case by Step 3a we know that 
$\partial(h^2_{k|R_k})$ is bounded from below near $x$.

We start with the case where $\partial h_k^1$ is bounded from below; in
other words, we are not dealing with tangency points but only with cusps.
In that case, we can use an argument similar to Step 3a, except that the
roles of the two components of $h_k$ are reversed. Namely, after a rotation
we assume that $\partial h_k^1(x)=\lambda dz_k^1$ for some nonzero constant
$\lambda$, and we define components $u_2,\dots,u_n$ of $\mathrm{Jac}(f_k)$
as previously (using $\partial h_k^1$ rather than $\partial h_k^2$ to define
the $\mu_k^j$). Let $\theta=\partial u_2\wedge\dots\wedge \partial u_n\wedge
\partial h_k^1$~: along $R_k$, the ratio between $\theta$ and
$\partial (h^1_{k|R_k})$, or equivalently $\partial(\phi_{k|R_k})$, is
bounded between two fixed constants, so the transverse vanishing of $\theta$
is what we are trying to obtain. More precisely, our aim is to show that
the quantity $(u_2,\dots,u_n,\theta)$, which is a section of a rank $n$
bundle $\mathcal{E}$ near $x$, can be made uniformly transverse to $0$ by
a small perturbation.

For this purpose, we first show the existence of complex-valued polynomials
$(P^1_j,P^2_j)$ and local sections $\epsilon_j$ of $\mathcal{E}$, 
$2\le j\le n+1$, such that~: 

(a) for any coefficients $w_j\in\C$, replacing the 
given sections of $L^{\otimes k}$ by $(s_k^0,s_k^1+\sum w_j P^1_j 
s_{k,x}^\mathrm{ref},s_k^2+\sum w_j P^2_j s_{k,x}^\mathrm{ref})$ affects 
$(u_2,\dots,u_n,\theta)$ by the addition of $\sum w_j \epsilon_j +
O(w_j^2)$~; 

(b) the sections $\epsilon_j$ define a local frame in
$\mathcal{E}$, and $\epsilon_2\wedge\dots\wedge\epsilon_{n+1}$ is bounded
from below by a universal constant.

By the same argument as in Step 3a, we find after a suitable rotation an
index $j_0\neq 1$ such that, letting $P^1_{n+1}=0$ and 
$P^2_{n+1}=z_k^2z_k^{j_0}$, the corresponding local section $\epsilon_{n+1}$
of $\mathcal{E}$ is, at $x$, of the form $(0,\dots,0,\epsilon'_{n+1})$, with
$\epsilon'_{n+1}$ bounded from below by a fixed constant.

Moreover, adding $w\,z_k^j s_{k,x}^\mathrm{ref}$ to $s_k^2$ amounts to
adding $w$ to $u_j$ and does not affect the other $u_i$'s, by the argument
in Step 2. So, letting $P^1_j=0$ and $P^2_j=z_k^j$, we get that the
corresponding local sections of $\mathcal{E}$ are of the form
$\epsilon_j=(0,\dots,1,\dots,0,\epsilon'_j)$, where the coefficient $1$ is
in $j$-th position.

So it is easy to check that both conditions (a) and (b) are satisfied by
these perturbations. The rest of the argument is as in Step 3a~:
expressing $(u_2,\dots,u_n,\theta)$ as a linear combination of $\epsilon_2,
\dots,\epsilon_{n+1}$, one uses Theorem 12 of \cite{D2} to obtain
transversality to $0$ over a small ball centered at $x$.

We now consider the second possibility, namely the case where $\partial
h_k^1(x)$ is small, which corresponds to tangency points. By property $(2)$
we know that $\partial h_k^2(x)$ is bounded from below, and we can assume
that it is colinear to $dz_k^1$. We then define components $u_2,\dots,u_n$
of $\mathrm{Jac}(f_k)$ as usual (as in Step 3a and unlike the previous case, 
the $\mu_k^j$ are defined using $\partial h_k^2$ rather than $\partial
h_k^1$). Letting $\theta=\partial u_2\wedge\dots\wedge \partial u_n\wedge
\partial h_k^1$, we want as before to obtain the transversality to $0$ of
the quantity $(u_2,\dots,u_n,\theta)$, which is a local section of a rank 
$n$ bundle $\mathcal{E}$ near $x$. For this purpose, as usual we look for
polynomials $P^1_j$, $P^2_j$ and local sections $\epsilon_j$ satisfying
the same properties (a) and (b) as above.

In order to construct $P^i_{n+1}$, observe that, by the result of Step 3a,
the quantity $\partial u_2\wedge\dots\wedge\partial u_n\wedge\partial h_k^2$
is bounded from below at $x$. So, adding to $s_k^1$ a small multiple of
$s_k^2$ does not affect the $u_j$'s, but it affects $\theta$ non-trivially.
However, this perturbation is not localized, so it is not suitable for our
purposes (we can't apply the globalization argument). Instead, let
$P^1_{n+1}$ be a polynomial of degree $2$ in the coordinates $z_k^j$ 
and their complex conjugates, such that $P^1_{n+1}s_{k,x}^\mathrm{ref}$ 
coincides with $s_k^2$ up to order two at $x$. Note that the coefficients of
$P^1_{n+1}$ are bounded by uniform constants, and that its antiholomorphic
part is at most of the order $O(k^{-1/2})$ (because $s_k^2$ and
$s_{k,x}^\mathrm{ref}$ are asymptotically holomorphic); therefore,
$P^1_{n+1}s_{k,x}^\mathrm{ref}$ is an admissible localized asymptotically
holomorphic perturbation. Also, define $P^2_{n+1}=0$.
Then one easily checks that the local section $\epsilon_{n+1}$ of 
$\mathcal{E}$ corresponding to $P^1_{n+1}$ and $P^2_{n+1}$ is, at $x$, of
the form $(0,\dots,0,\epsilon'_{n+1})$, where $\epsilon'_{n+1}$ is bounded
from below. 

Moreover, let $P^1_j=z_k^j$ and $P^2_j=0$~: as above, this perturbation
affects $u_j$ and not the other $u_i$'s, and we get that the
corresponding local sections of $\mathcal{E}$ are of the form
$\epsilon_j=(0,\dots,1,\dots,0,\epsilon'_j)$, where the coefficient $1$ is
in $j$-th position.

Once again, these perturbations satisfy both conditions (a) and (b).
Therefore, expressing $(u_2,\dots,u_n,\theta)$ as a linear combination of 
$\epsilon_2,\dots,\epsilon_{n+1}$, Theorem 12 of \cite{D2} yields
transversality to $0$ over a small ball centered at $x$ by the usual
argument. Now that both possible cases have been handled, we can apply the
standard globalization argument to obtain uniform transversality to the
stratum $S\subset\mathcal{J}^2(X,\CP^2)$. This gives properties $(4)$ and
$(6)$ of Definition 3.2.

{\bf Step 4.} Now that all required transversality properties have been
obtained, we perform further perturbations in order to achieve the other
conditions in Definition 3.2. These new perturbations are bounded by a 
fixed multiple of $k^{-1/2}$, so the transversality properties are not
affected. The argument is almost the same as in the case of 4-manifolds
(see \S 4 of \cite{A2} and \S 3.1 of \cite{AK}); the adaptation to the 
higher-dimensional case is very easy.

One first defines a suitable almost-complex structure $\tilde{J}_k$, by the
same argument as in \S 4.1 of \cite{A2} (except that one also considers
the points of $\mathcal{T}_k$ and $\mathcal{I}_k$ besides the cusps).
As explained in \S 4.1 of \cite{A2}, a suitable perturbation makes 
it possible to obtain the local holomorphicity of $f_k$ near these points, 
which yields conditions $(4')$, $(6')$ and $(8')$~;
the argument is the same in all three cases. Next, a
generically chosen small perturbation yields the self-transversality of $D$ 
(property $(7)$). Finally, as
described in \S 4.2 of \cite{A2}, a suitable perturbation yields property 
$(3')$ along the branch curve without modifying $R_k$ and $D_k$ and
without affecting the other compatibility properties. This completes the
proof of the existence statement in Theorem 4.1.

{\bf Uniqueness.}
The uniqueness statement is obtained by showing that, provided that $k$ is 
large enough, the whole argument extends to the case of families of sections 
depending continuously on a parameter $t\in [0,1]$. Then, given two sequences 
of quasiholomorphic maps, one can start with a one-parameter family of 
sections interpolating between them in a trivial way and perturb it in such 
a way that the required properties hold for all parameter values (with the 
exception of $(7)$ when a node cancellation occurs). If one moreover checks
that the construction can be performed in such a way that the two end points 
of the one-parameter family are not affected by the perturbation, the 
isotopy result becomes an immediate corollary. Observe that, in the
one-parameter construction, the almost-complex structure is allowed to 
depend on $t$.

Most of the above argument extends to 1-parameter families in a 
straightforward manner, exactly as in the four-dimensional case~; the key
observation is that all the standard building blocks (existence
of approximately holomorphic Darboux coordinates $z_k^j$ and of localized
approximately holomorphic sections $s_{k,x}^\mathrm{ref}$, local
transversality result, globalization principle, ...) remain valid in the
parametric case, even when the almost-complex structure depends on $t$.
The only places where the argument differs from the case of 4-manifolds are 
properties $(3)$, $(4)$ and $(6)$, obtained in Steps 2 and 3 above. 

For property $(3)$, one easily checks that it is still possible in the
parametric case to assume, after composing with suitable rotations depending
continuously on the parameter $t$, that $s_k^1(x)=s_k^2(x)=0$ and that
$\partial h_k^1(x)$ is bounded from below and directed along $dz_k^1$. This
makes it possible to define $\mu_k^j$ and $u_k^j$ as in the non-parametric
case, and the parametric version of Theorem 12 of \cite{D2} yields a
suitable perturbation depending continuously on $t$.

The argument of Step 3a also extends to the parametric case, using the
following observation. Fix a point $x\in X$, and let
$\rho_k(t)=|\partial\phi_{k,t}(x)|$. For all values of $t$ such that 
$\rho_k(t)$ is small enough (smaller than a fixed constant $\alpha>0$), 
we can perform the construction as in the non-parametric case, defining 
$u_{j,t}$ and $\theta_t$. If $\rho'_k(t)=|\theta_t(x)|$ is small enough 
(smaller than $\alpha$), then we can apply the same argument as in the 
non-parametric case to define polynomials $(P^1_{j,t},P^2_{j,t})$ and local
sections $\epsilon_{j,t}$ of $\mathcal{E}_0$. However the definition of
$P^1_{n+1}$ needs to be modified as follows. Although it
is still possible after a suitable rotation depending continuously on $t$
to assume that $\partial u_2\wedge\partial h_k^2(x)$ is small, the 
choice of an index $j_0\neq 1$ such that $dz_k^{j_0}$ lies far from the span 
of the $\partial u_j(x)$ may depend on $t$. Instead, we define $\nu_{k,t}$ 
as a unit vector in $\C^{n-1}$ depending continuously on $t$ and such that
$\sum_{j=2}^n \nu_{k,t}^j\,dz_k^j$ lies far from the span of $\partial
u_j(x)$, and let $P^1_{n+1,t}=\sum_{j=2}^n \nu_{k,t}^j z_k^2 z_k^j$. Then
the required properties are satisfied, and we can proceed with the argument.
So, provided that $\rho_k(t)$ and $\rho'_k(t)$ are both smaller than
$\alpha$, we can use Theorem 12 of \cite{D2} to obtain a localized 
perturbation $\tau_{k,t}$ depending continuously on $t$ and such that 
$s_{k,t}+\tau_{k,t}$ satisfies the desired transversality property near $x$.

In order to obtain a well-defined perturbation for all values of $t$, we
introduce a continuous cut-off function $\beta:\R_+\to [0,1]$ which equals
$1$ over $[0,\alpha/2]$ and vanishes outside of $[0,\alpha]$. Then, we set
$\tilde\tau_{k,t}=\beta(\rho_k(t))\beta(\rho'_k(t))\tau_{k,t}$, which is
well-defined for all $t$ and depends continuously on $t$. Since
$s_{k,t}+\tilde\tau_{k,t}$ coincides with $s_{k,t}+\tau_{k,t}$ when
$\rho_k(t)$ and $\rho'_k(t)$ are smaller than $\alpha/2$, the
required transversality holds for these values of $t$~; moreover, for the
other values of $t$ we know that the 2-jet of $s_{k,t}$ already lies at
distance more than $\alpha/2$ from the stratum $S_0$, and we can safely
assume that $\tilde\tau_{k,t}$ is much smaller than $\alpha/2$, so the
perturbation does not affect transversality.
Therefore we obtain a well-defined local perturbation for all $t\in [0,1]$,
and the one-parameter version of the result of Step 3a follows by the
standard globalization argument.

The argument of Step 3b is extended to one-parameter families in the same
way~: given a point $x\in X$, the same ideas as for Step 3a yield, for 
all values of the parameter $t$ such that the 2-jet of $s_{k,t}$ at $x$
lies close to the stratum $S$, small localized perturbations $\tau_{k,t}$
depending continuously on $t$ and such that $s_{k,t}+\tau_{k,t}$ satisfies
the desired property over a small ball centered at $x$. As seen above, two
different types of formulas for $\tau_{k,t}$ arise depending on which
component of the stratum $S$ is being hit; however, the result of Step 3a 
implies that, in any interval of parameter values such that the jet of 
$s_{k,t}$ remains close to $S$, only one of the two components of $S$ has to 
be considered, so $\tau_{k,t}$ indeed depends continuously on $t$. The same
type of cut-off argument as for Step 3a then makes it possible to extend
the definition of $\tau_{k,t}$ to all parameter values and complete the 
proof.
\end{proof}

\subsection{The topology of quasiholomorphic maps}

We now describe the topological features of quasiholomorphic maps and the 
local models which characterize them near the critical points.

\begin{prop}
Let $f_k:X-Z_k\to\CP^2$ be a sequence of quasiholomorphic maps. Then the
fibers of $f_k$ are codimension $4$ symplectic submanifolds, intersecting
at the set of base points $Z_k$, and smooth away from the critical curve
$R_k\subset X$. The submanifolds $R_k$ and $Z_k$ of $X$ are smooth and 
symplectic, and the image $f_k(R_k)=D_k$ is a symplectic braided curve
in $\CP^2$. 

Moreover, given any point $x\in R_k$, there exist local approximately
holomorphic coordinates on $X$ near $x$ and on $\CP^2$ near $f_k(x)$ in
which $f_k$ is topologically conjugate to one of the two following models~:

$(i)$ $(z_1,\dots,z_n)\mapsto (z_1^2+\dots+z_{n-1}^2,z_n)$
$($points where $f_{k|R_k}$ is an immersion$)$~;

$(ii)$ $(z_1,\dots,z_n)\mapsto (z_1^3+z_1z_n+z_2^2+\dots+z_{n-1}^2,z_n)$
$($near the cusp points$)$.
\end{prop}

\begin{proof} The smoothness and symplecticity properties of the various
submanifolds appearing in the statement follow from the observation made
by Donaldson in \cite{D1} that the zero sets of approximately holomorphic
sections satisfying a uniform transversality property are smooth and
approximately $J$-holomorphic, and therefore symplectic. In particular,
the smoothness and symplecticity of the fibers of $f_k$ away from $R_k$
follow immediately from Definition 3.2~: since $\mathrm{Jac}(f_k)$ is
bounded from below away from $R_k$ (because it satisfies a uniform
transversality property), and since the sections $s_k$ are asymptotically
holomorphic, it is easy to check that the level sets of $f_k$ are, away from
$R_k$, smooth symplectic submanifolds. Symplecticity near the singular points 
is an immediate consequence of the local models $(i)$ and $(ii)$ that we will 
obtain later in the proof.

The corresponding properties of $Z_k$ and $R_k$ are obtained by the
same argument~: $Z_k$ and $R_k$ are the zero sets of asymptotically
holomorphic sections, both satisfying a uniform transversality property
(by conditions $(1)$ and $(3)$ of Definition 3.2, respectively), so they
are smooth and symplectic.\medskip

We now study the local models at critical points of $f_k$. We start with
the case of a cusp point $x\in X$. By property $(2)$ of Definition 3.2, 
$\partial f_k$ has complex rank $1$ at $x$, so we can find local complex 
coordinates $(Z_1,Z_2)$ on $\CP^2$ near $f_k(x)$ such that 
$\mathrm{Im}\,\partial f_k(x)$ is the $Z_2$ axis. Pulling back $Z_2$ via 
the map $f_k$, we obtain, using property $(4')$, a
$\tilde{J}_k$-holomorphic function whose differential does not vanish 
near $x$~; therefore, we can find a $\tilde{J}_k$-holomorphic 
coordinate chart $(z_1,\dots,z_n)$ on $X$ at $x$ such that $z_n=Z_2\circ f_k$.
In the chosen coordinates, we get $f_k(z_1,\dots,z_n)=(g(z_1,\dots,z_n),z_n)$,
where $g$ is holomorphic and $\partial g(0)=0$.

Since $x$ is by assumption a cusp point, the tangent direction to 
$R_k$ at $x$ lies in the kernel of $\partial f_k(0)$, i.e.\ in the span of
the $n-1$ first coordinate axes~; after a suitable rotation we may assume
that $T_xR_k$ is the $z_1$ axis.
Near the origin, $\mathrm{Jac}(f_k)$ is characterized by its $n-1$
components $(\partial g/\partial z_1,\dots,\partial g/\partial z_{n-1})$,
and the critical curve $R_k$ is the set of points where these quantities
vanish. Therefore, at the origin, 
$\partial^2 g/\partial z_1^2 =\partial^2 g/\partial z_1\partial z_2=\dots=
\partial^2 g/\partial z_1\partial z_{n-1}=0$.
Nevertheless, $\mathrm{Jac}(f_k)$ vanishes transversely to $0$ at the
origin, so the matrix of second derivatives $M=(\partial^2 g/
\partial z_i\partial z_j(0))$, $2\le i\le n$, $1\le j\le n-1$, is
non-degenerate (invertible) at the origin. In particular, the first column
of $M$ (corresponding to $j=1$) is non-zero, and therefore
$\partial^2 g/\partial z_1\partial z_n(0)$ is necessarily non-zero~; after a
suitable rescaling of the coordinates we may assume that this coefficient is
equal to $1$. Moreover, the invertibility of $M$ implies that the submatrix
$M'=(\partial^2 g/\partial z_i\partial z_j(0))$, $2\le i,j\le n-1$ is 
also invertible, i.e.\ it represents a non-degenerate quadratic form.

Diagonalizing this quadratic form, we can assume after a suitable linear
change of coordinates that the diagonal coefficients of $M'$ are equal to
$2$ and the others are zero. Therefore $g$ is of the form
$g(z_1,\dots,z_n)=z_1z_n+\sum_{j=2}^{n-1} z_j^2+\sum_{j=2}^{n-1}
\alpha_jz_jz_n+O(z^3)$. Changing coordinates on $X$ to replace $z_j$ by
$z_j+\frac{1}{2}\alpha_jz_n$ for all $2\le j\le n-1$, and on $\CP^2$ to
replace $Z_1$ by $Z_1+\frac{1}{4}\sum \alpha_j^2 Z_2^2$, we can ensure
that $g(z_1,\dots,z_n)=z_1z_n+\sum_{j=2}^{n-1} z_j^2+O(z^3)$.

Observe that $R_k$ is described near the origin by expressing the
coordinates $z_2,\dots,z_n$ as functions of $z_1$. By assumption the 
expressions of $z_2,\dots,z_n$ are all of the form $O(z_1^2)$. Substituting
into the formula for $\mathrm{Jac}(f_k)$, and letting $g_{ijk}=\partial^3g/
\partial z_i\partial z_j\partial z_k(0)$, we get that local equations of 
$R_k$ near the origin are $z_j=-\frac{3}{2}g_{j11}z_1^2+O(z_1^3)$ for
$2\le j\le n-1$, and $z_n=-3g_{111}z_1^2+O(z_1^3)$. It follows that
$f_{k|R_k}$ is locally given in terms of $z_1$ by the map
$z_1\mapsto (-2g_{111}z_1^3+O(z_1^4),-3g_{111}z_1^2+O(z_1^3))$.
Therefore, the transverse vanishing of $\partial(f_{k|R_k})$ at the
origin implies that $g_{111}\neq 0$, so after a suitable rescaling we may
assume that the coefficient of $z_1^3$ in the power series expansion of $g$
is equal to one.

On the other hand, suitable coordinate changes can be used to kill all
other degree~$3$ terms in the expansion of $g$~: if $2\le i\le
n-1$ the coefficient of $z_iz_jz_k$ can be made zero by replacing $z_i$ by 
$z_i+\frac{c}{2}z_jz_k$~; similarly for $z_n^3$ (replace $Z_1$ by
$Z_1+cZ_2^3$), $z_1z_n^2$ and $z_1^2z_n$ (replace $z_1$ by 
$z_1+cz_n^2+c'z_1z_n$). So we get that $f_k(z_1,\dots,z_n)=(z_1^3+z_1z_n+
z_2^2+\dots+z_{n-1}^2+O(z^4),z_n)$. It is then a standard result of
singularity theory that the higher order terms can be absorbed by suitable
coordinate changes.
\medskip

We now turn to the case of where $x$ is a point of $R_k$ which does not lie
close to any of the cusp points. Conditions $(2)$ and $(3')$ imply that the 
differential of $f_k$ at $x$ 
has real rank $2$ and that its image lies close to a complex line in the
tangent plane to $\CP^2$ at $f_k(x)$. Therefore, there exist local
approximately holomorphic coordinates $(Z_1,Z_2)$ on $\CP^2$ such that 
$\mathrm{Im}\,\nabla f_k(x)$ is the $Z_2$ axis. Moreover, because 
$Z_2\circ f_k$ is an approximately holomorphic function whose derivative
at $x$ satisfies a uniform lower bound, it remains possible to find local 
approximately holomorphic coordinates $z_1,\dots,z_n$ on $X$ such that 
$z_n=Z_2\circ f_k$. As before, we can write
$f_k(z_1,\dots,z_n)=(g(z_1,\dots,z_n),z_n)$,
where $g$ is an approximately holomorphic function such that $\nabla g(0)=0$.

By assumption $f_k$ restricts to $R_k$ as
an immersion at $x$, so the projection to the $z_n$ axis of $T_xR_k$ is
non-trivial. In fact, property $(4)$ implies that, if $\partial(f_{k|R_k})$
is very small at $x$, then a cusp point lies nearby~; so we can assume that
the $z_n$ component of $T_xR_k$ is larger than some fixed constant. As a
consequence, one can show that $R_k$ is locally given by equations of the 
form $z_j=h_j(z_n)$, where the functions $h_j$ are approximately holomorphic 
and have bounded derivatives. Therefore, a suitable change of coordinates on
$X$ makes it possible to assume that $R_k$ is locally given by the equations
$z_1=\dots=z_{n-1}=0$. Similarly, a suitable approximately holomorphic
change of coordinates on $\CP^2$ makes it possible to assume that $f_k(R_k)$ 
is locally given by the equation $Z_1=0$.

As a consequence, we have that $g_{|R_k}=0$ and, since the image of $\nabla
f_k$ at a point of $R_k$ coincides with the tangent space to $f_k(R_k)$,
$\nabla g$ vanishes at all points of $R_k$. In particular this implies that
$\partial^2 g/\partial z_j\partial z_n(0)=0$ for all $1\le j\le n$.
Moreover, property $(3)$ implies that $\mathrm{Jac}(f_k)$ 
vanishes transversely at the origin, and therefore that the matrix 
$(\partial^2 g/\partial z_i\partial z_j(0))$, $1\le i,j\le n-1$
is invertible, i.e.\ it represents a non-degenerate quadratic form.
This quadratic form can be diagonalized by a suitable change of
coordinates~; because the transversality property $(3)$ is uniform, the 
coefficients are bounded between fixed constants. After a suitable
rescaling, we can therefore assume that $\partial^2 g/\partial z_i
\partial z_j(0)$ is equal to $2$ if $i=j$ and $0$ otherwise. 

In conclusion, we get that
$g(z_1,\dots,z_n)=z_1^2+\dots+z_{n-1}^2+h(z_1,\dots,z_n)$, where $h$
is the sum of a holomorphic function which vanishes up to order $3$ at 
the origin and of a non-holomorphic function which vanishes up to order
$2$ at the origin and has derivatives bounded by $O(k^{-1/2})$.

Let $z$ be the column vector $(z_1,\dots,z_{n-1})$, and denote by
$\mathbf{z}$ the vector $(z_1,\dots,z_n)$. Using the fact that $g$ 
vanishes up to order $2$ along $R_k$, we conclude that
there exist matrix-valued functions $\alpha$, $\beta$ and $\gamma$
with the following properties~:
\smallskip

$(a)$ $g(\mathbf{z})={}^t\!z\alpha(\mathbf{z}) z+{}^t\!\bar{z}
\beta(\mathbf{z}) z+ {}^t\!\bar{z}\gamma(\mathbf{z})\bar{z}$~; 
($\alpha$ and $\gamma$ are symmetric)~;

$(b)$ $\alpha$ is approximately holomorphic and has uniformly bounded 
derivatives~; $\alpha(0)=I$~;

$(c)$ $\beta$ and $\gamma$ and their derivatives are bounded by fixed
multiples of $k^{-1/2}$.\smallskip

\noindent
The implicit function theorem then makes it possible to 
construct a $C^\infty$ approximately holomorphic change of coordinates 
of the form $z\mapsto\lambda(\mathbf{z})z+\mu(\mathbf{z})\bar{z}$ 
(with $\lambda(0)$ orthogonal, $\lambda$ approximately holomorphic, 
$\mu=O(k^{-1/2})$), such that $g$ becomes of the form $g(\mathbf{z})=
{}^t\!zz+{}^t\!\bar{z}\tilde\gamma(\mathbf{z})\bar{z}$.

Unfortunately, smooth coordinate changes are not sufficient to further
simplify this expression; instead, in order to obtain the desired local
model one must use as coordinate change an ``approximately holomorphic 
homeomorphism'', which is smooth away from $R_k$ but admits only 
directional derivatives at the points of $R_k$. More precisely,
starting from $g={}^t\!zz+h$ and using that $h/|z|^2$ is bounded by
$O(k^{-1/2})+O(\mathbf{z})$, we can write
$$g(\mathbf{z})=\sum_{j=1}^{n-1}\tilde{z}_j^2,\qquad \tilde{z}_j=z_j
\Bigl(1+\frac{\bar{z}_j}{z_j}\frac{h(\mathbf{z})}{|z|^2}\Bigr)^{1/2}.$$
This gives the desired local model and ends the proof.
\end{proof}

\begin{rem}
The local model at points of $R_k$ only holds topologically (up to
an approximately holomorphic homeomorphism), which is not fully satisfactory.
However, by replacing $(3')$ by a stronger condition, it is possible to 
obtain the same result in smooth approximately holomorphic coordinates.
This new condition can be formulated as follows. Away from the cusp points,
the complex lines $(\mathrm{Im}\,\partial f_k)^\bot$ define a line bundle 
$V\subset T\CP^2_{|D_k}$, everywhere transverse to $TD_k$. A neighborhood of 
the zero section in $V$ can be sent via the exponential map of the
Fubini-Study metric onto a neighborhood of $D_k$ (away from the cusps), in 
such a way that each fiber $V_x$ is mapped holomorphically to a subset
$\mathcal{V}_x$ contained in a complex line in $\CP^2$. 

Lifting back to a neighborhood of $R_k$ in $X$,
we can define slices $\mathcal{W}_x=f_k^{-1}(\mathcal{V}_{f_k(x)})$ for all
$x\in R_k$ lying away from $\mathcal{C}_k$. It is then possible to identify
a neighborhood of $R_k$ (away from $\mathcal{C}_k$) with a neighborhood of
the zero section in the vector bundle $W$ whose fiber at $x\in R_k$ is
$\mathrm{Ker}\,\partial f_k(x)$, in such a way that each fiber $W_x$ gets mapped to
$\mathcal{W}_x$. Observe moreover that, since $W_x$ is a complex subspace
in $(T_xX,\tilde{J}_k)$, $W$ is endowed with a natural complex
structure induced by $\tilde{J}_k$. It is then possible to ensure that the
``exponential map'' from $W_x$ to $\mathcal{W}_x$ is approximately
$\tilde{J}_k$-holomorphic for every $x$, and, using condition $(4')$, 
holomorphic when $x$ lies at distance less than $\delta/2$ from a cusp
point.

With this setup understood, and composing on both sides with the exponential
maps, $f_k$ induces a fiber-preserving map $\psi_k$ between the bundles $W$ 
and $V$~; this map is approximately holomorphic everywhere, and holomorphic 
at distance less than $\delta/2$ from $\mathcal{C}_k$. The condition which
we impose as a replacement of $(3')$ is that $\psi_k$ should be fiberwise
holomorphic over a neighborhood of the zero section in $W$. 

The proof of existence of quasiholomorphic maps satisfying this strengthened
condition follows a standard argument~: trivializing locally $V$ and $W$ for
each value of $k$, and given asymptotically holomorphic maps $\psi_k$, 
Lemma 8 of \cite{A2} (see also \cite{D1}) implies the existence of a
fiberwise holomorphic map $\tilde\psi_k$ differing from $\psi_k$ by
$O(k^{-1/2})$ over a neighborhood of the zero section. It is moreover
easy to check that $\tilde\psi_k=\psi_k$ near the cusp points. So, in order
to obtained the desired property, we introduce a smooth cut-off function and
define a map $\hat{\psi}_k$ which equals $\tilde\psi_k$ near the zero
section and coincides with $\psi_k$ beyond a certain distance. 
Going back through the exponential maps, we obtain a map $\hat{f}_k$ which
differs from $f_k$ by $O(k^{-1/2})$ and coincides with $f_k$ outside a small
neighborhood of $R_k$ and near the cusp points. The corresponding perturbations of the 
asymptotically holomorphic sections $s_k\in\Gamma(\C^3\otimes L^{\otimes k})$
are easy to construct. Moreover, we can always assume that $\tilde\psi_k$ and 
$\psi_k$ coincide at order $1$ along the zero section, i.e.\ that
$\hat{f}_k$ and $f_k$ coincide up to order 1 along the branch curve~;
therefore, the branch curve of $\hat{f}_k$ and its image are the same as for
$f_k$, and so all properties of Definition 3.2 hold for $\hat{f}_k$.

Once this condition is satisfied, getting the correct local model at a point
$x\in R_k$ in smooth approximately holomorphic coordinates is an easy task.
Namely, we can define, near $f_k(x)$, local approximately holomorphic 
coordinates $Z_2$ on $D_k$ and $Z_1$ on the fibers of $V$ ($Z_1$ is a
complex linear function on each fiber, depending approximately
holomorphically on $Z_2$). Using the exponential map, we can use $(Z_1,Z_2)$ 
as local coordinates on $\CP^2$. Lifting $Z_2$ via
$\hat{f}_k$ yields a local coordinate $z_n$ on $R_k$ near $x$. Moreover, we
can locally define complex linear coordinates $z_1,\dots,z_{n-1}$ in the 
fibers of $W$, depending approximately holomorphically on $z_n$. Using again
the exponential map, $(z_1,\dots,z_n)$ define local approximately
holomorphic coordinates on $X$. Then, by construction, local equations are
$z_1=\dots=z_{n-1}=0$ for $R_k$ and $Z_1=0$ for $D_k$, and $f_k$ is given by
$f_k(z_1,\dots,z_n)=(\psi_k(z_1,\dots,z_n),z_n)$. Moreover, we know that
$\psi_k$ is, for each value of $z_n$, a holomorphic function of
$z_1,\dots,z_{n-1}$, vanishing up to order 2 at the origin. We can then
use the argument in the proof of Proposition 4.2 to obtain the expected
local model in smooth approximately holomorphic coordinates. 
\end{rem}

\subsection{Monodromy invariants of quasiholomorphic maps}

We now look at the monodromy invariants naturally arising from
quasiholomorphic maps to $\CP^2$. Let $f:X-Z\to\CP^2$ be one of the maps
constructed in Theorem 3.1 for large enough $k$. The fibers of $f$ are 
singular along the smooth symplectic curve $R\subset X$, whose image in 
$\CP^2$ is a symplectic braided curve. Therefore, we obtain a first 
interesting invariant by considering the critical curve $D\subset \CP^2$.

As in the four-dimensional case, using the projection $\pi:\CP^2-\{(0\!:\!0
\!:\!1)\}\to\CP^1$ we can describe the topology of $D$ by 
a braid monodromy map \begin{equation}
\rho_n:\pi_1(\C-\{p_1,\dots,p_r\})\to B_d,\end{equation}
where $p_1,\dots,p_r$ are
the images by $\pi$ of the cusps, nodes and tangency points of 
$D$, and $d=\deg D$. Alternately, we can also express this monodromy
as a braid group factorization
\begin{equation}\Delta^2=\prod_{j=1}^r Q_j^{-1}X_1^{r_j}Q_j.\end{equation}
Like in the four-dimensional case, this braid factorization completely
characterizes the curve $D$ up to isotopy, but it is only well-defined
up to simultaneous conjugation and Hurwitz equivalence. 

We now turn to the second part of the problem, namely describing the
topology of the map $f:X-Z\to\CP^2$ itself. As in the case of Lefschetz
pencils, we blow up $X$ along $Z$ in order to obtain a well-defined map
$\hat{f}:\hat{X}\to\CP^2$. The fibers of $\hat{f}$ are naturally identified
with those of $f$, made mutually disjoint by the blow-up process.

Denote by $\Sigma^{2n-4}$ the generic fiber, i.e.\ the fiber 
above a point of $\CP^2-D$. The structure of the singular
fibers of $\hat{f}$ can be easily understood by looking at the local models
obtained in Proposition 4.2. The easiest case is that of the fiber above a
smooth point of $D$. This fiber intersects $R$ transversely in one
point, where the local model is $(z_1,\dots,z_n)\mapsto
(z_1^2+\dots+z_{n-1}^2,z_n)$, which can be thought of as a
one-parameter version of the model map for the singularities of a Lefschetz 
pencil in dimension $2n-2$. Therefore, as in that case, the singular fiber 
is obtained by collapsing a vanishing cycle, namely a Lagrangian sphere 
$S^{n-2}$, in the generic fiber $\Sigma$, and
the monodromy of $\hat{f}$ maps a small loop around 
$D$ to a positive Dehn twist along the vanishing cycle.

The fiber of $\hat{f}$ above a nodal point of $D$ intersects $R$
transversely in two points, and is similarly obtained from $\Sigma$ by 
collapsing two disjoint Lagrangian spheres. In fact, the nodal point does
not give rise to any specific local model in $X$, as it simply corresponds 
to the situation where two points of $R$ happen to lie in the same fiber. 

Finally, in the case of a cusp
point of $D$, the local model $(z_1,\dots,z_n)\mapsto
(z_1^3+z_1z_n+z_2^2+\dots+z_{n-1}^2,z_n)$ can be used to show that the
singular fiber is a ``fishtail'' fiber, obtained by collapsing two Lagrangian 
spheres which intersect transversely in one point.

With this understood, the topology of $\hat{f}$ is described by its monodromy 
around the
singular fibers. As in the case of Lefschetz fibrations, the monodromy 
consists of symplectic automorphisms of $\Sigma$ preserving the submanifold 
$Z$. However, as in \S 2, defining a monodromy map with values in
$\mathrm{Map}^\omega(\Sigma,Z)$ requires a trivialization of the normal
bundle of $Z$, which is only possible over an affine subset
$\C^2\subset\CP^2$. So, the monodromy of $\hat{f}$ is described by a group
homomorphism
\begin{equation}\psi_n:\pi_1(\C^2-D)\to\mathrm{Map}^\omega(\Sigma,Z).
\end{equation}

A simpler description can be obtained by restricting oneself to a generic
line $L\subset\CP^2$ which intersects $D$ transversely in $d$ points
$q_1,\dots,q_d$. In fact, Definition 3.2 implies that we 
can use the fiber of $\pi$ above $(0\!:\!1)$ for this purpose.
As in \S 3.2, the inclusion $i:\C-\{q_1,\dots,q_d\}\to\C^2-D$ induces
a surjective homomorphism on fundamental groups. The relations between
the geometric generators $\gamma_1,\dots,\gamma_d$ of $\pi_1(\C^2-D)$ are
again given by the braid factorization (one relation for each factor) in the
same manner as in \S 3.2. Note that the relation 
$\gamma_1\dots\gamma_d=1$ only holds in $\pi_1(\CP^2-D)$, not in 
$\pi_1(\C^2-D)$.

It follows from these observations that the monodromy of $\hat{f}$ can be
described by the monodromy morphism
\begin{equation}
\theta_{n-1}:\pi_1(\C-\{q_1,\dots,q_d\})\to
\mathrm{Map}^\omega(\Sigma,Z)
\end{equation}
defined by $\theta_{n-1}=\psi_n\circ i_*$. We know from the above discussion 
on the structure of $\hat{f}$ near its critical points that $\theta_{n-1}$
maps the geometric generators of $\pi_1(\C-\{q_1,\dots,q_d\})$ to positive
Dehn twists. Moreover, by considering the normal bundle to the
exceptional divisor in $\hat{X}$ one easily checks that the monodromy
around infinity is again a twist along $Z$ in $\Sigma$, 
i.e.\ $\theta_{n-1}(\gamma_1\dots\gamma_d)=\delta_Z$. 

These properties of $\theta_{n-1}$ are strikingly similar to those of the 
monodromy of a symplectic Lefschetz pencil. In fact, 
let $W=f^{-1}(L)$ be the preimage of a complex line 
$L=\CP^1\subset\CP^2$ intersecting $D$ transversely. Then
the restriction of $f$ to the smooth symplectic hypersurface $W\subset X$ 
endows it with a structure of
symplectic Lefschetz pencil with generic fiber $\Sigma$ and base set $Z$~; 
for example, if one chooses $L=\pi^{-1}(0\!:\!1)$, then $W$ is the zero set 
of $s_k^0$ and the restricted pencil $f_{|W}:W-Z\to\CP^1$ is defined by the 
two sections $s_k^1$ and $s_k^2$. The monodromy of the restricted pencil 
is, by construction, given by the map $\theta_{n-1}$.

The situation is summarized in the following picture~:

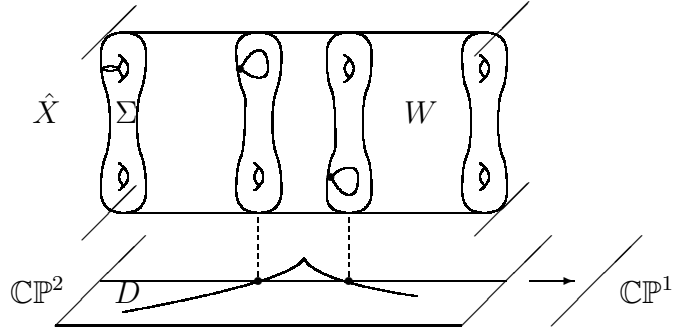
\setlength{\unitlength}{0.6mm}
\begin{center}
\begin{picture}(120,78)(-50,-28)

\put(-40,0){\line(1,0){80}}
\put(-40,40){\line(1,0){80}}
\put(-45,-15){\line(1,0){90}}
\put(-65,-20){$\CP^2$}
\put(-55,-25){\line(1,0){90}}
\put(-55,-25){\line(1,1){20}}
\put(35,-25){\line(1,1){20}}
\put(55,-25){\line(1,1){20}}
\put(-49,34){\line(1,1){12}}
\put(-49,-6){\line(1,1){12}}
\put(38,-5.2){\line(1,1){12}}
\put(38,34.8){\line(1,1){12}}
\put(-60,20){$\hat{X}$}
\put(70,-20){$\CP^1$}
\put(50,-15){\vector(1,0){10}}

\qbezier(-40,40)(-35,40)(-35,35)
\qbezier(-35,35)(-35,30)(-36,27.5)
\qbezier(-36,27.5)(-37,25)(-37,20)
\qbezier(-40,0)(-35,0)(-35,5)
\qbezier(-35,5)(-35,10)(-36,12.5)
\qbezier(-36,12.5)(-37,15)(-37,20)
\qbezier(-40,40)(-45,40)(-45,35)
\qbezier(-45,35)(-45,30)(-44,27.5)
\qbezier(-44,27.5)(-43,25)(-43,20)
\qbezier(-40,0)(-45,0)(-45,5)
\qbezier(-45,5)(-45,10)(-44,12.5)
\qbezier(-44,12.5)(-43,15)(-43,20)
\qbezier(-41,35)(-37,32)(-41,29)
\qbezier(-40,34)(-42,32)(-40,30)
\qbezier(-41,5)(-37,8)(-41,11)
\qbezier(-40,6)(-42,8)(-40,10)

\qbezier(40,40)(35,40)(35,35)
\qbezier(35,35)(35,30)(36,27.5)
\qbezier(36,27.5)(37,25)(37,20)
\qbezier(40,0)(35,0)(35,5)
\qbezier(35,5)(35,10)(36,12.5)
\qbezier(36,12.5)(37,15)(37,20)
\qbezier(40,40)(45,40)(45,35)
\qbezier(45,35)(45,30)(44,27.5)
\qbezier(44,27.5)(43,25)(43,20)
\qbezier(40,0)(45,0)(45,5)
\qbezier(45,5)(45,10)(44,12.5)
\qbezier(44,12.5)(43,15)(43,20)
\qbezier(39,35)(43,32)(39,29)
\qbezier(40,34)(38,32)(40,30)
\qbezier(39,5)(43,8)(39,11)
\qbezier(40,6)(38,8)(40,10)

\qbezier(-10,40)(-5,40)(-5,35)
\qbezier(-5,35)(-5,30)(-6,27.5)
\qbezier(-6,27.5)(-7,25)(-7,20)
\qbezier(-10,0)(-5,0)(-5,5)
\qbezier(-5,5)(-5,10)(-6,12.5)
\qbezier(-6,12.5)(-7,15)(-7,20)
\qbezier(-10,40)(-15,40)(-15,36)
\qbezier(-15,36)(-15,32)(-14,28)
\qbezier(-14,28)(-13,26)(-13,25)
\qbezier(-10,0)(-15,0)(-15,5)
\qbezier(-15,5)(-15,10)(-14,12.5)
\qbezier(-14,12.5)(-13,15)(-13,25)
\qbezier(-8,32)(-8,36)(-10,36)
\qbezier(-10,36)(-12,36)(-14,32)
\qbezier(-8,32)(-8,30)(-10,30)
\qbezier(-10,30)(-12,30)(-14,32)
\qbezier(-11,5)(-7,8)(-11,11)
\qbezier(-10,6)(-12,8)(-10,10)

\qbezier(10,40)(15,40)(15,35)
\qbezier(15,35)(15,30)(14,27.5)
\qbezier(14,27.5)(13,25)(13,20)
\qbezier(10,0)(15,0)(15,5)
\qbezier(15,5)(15,10)(14,12.5)
\qbezier(14,12.5)(13,15)(13,20)
\qbezier(10,0)(5,0)(5,4)
\qbezier(5,4)(5,8)(6,12)
\qbezier(6,12)(7,14)(7,15)
\qbezier(10,40)(5,40)(5,35)
\qbezier(5,35)(5,30)(6,27.5)
\qbezier(6,27.5)(7,25)(7,15)
\qbezier(12,8)(12,4)(10,4)
\qbezier(10,4)(8,4)(6,8)
\qbezier(12,8)(12,10)(10,10)
\qbezier(10,10)(8,10)(6,8)
\qbezier(9,35)(13,32)(9,29)
\qbezier(10,34)(8,32)(10,30)

\put(20,20){ $W$}
\put(-44,20){ $\Sigma$}

\put(-42,-20){$D$}
\put(-14,32){\circle*{2}}
\put(6,8){\circle*{2}}
\qbezier(-44.8,32)(-42.9,31)(-41,32)
\qbezier[8](-44.8,32)(-42.9,33)(-41,32)
\multiput(-10,-1)(0,-2){7}{\line(0,-1){1}}
\multiput(10,-1)(0,-2){7}{\line(0,-1){1}}
\put(10,-15){\circle*{2}}
\put(-10,-15){\circle*{2}}
\qbezier(-40,-22)(-2,-15)(0,-10)
\qbezier(25,-18.5)(2,-15)(0,-10)
\end{picture}
\end{center}

\begin{rem}
If a cusp point of $D$ happens to lie close to the chosen line $L$, then 
two singular points of the restricted pencil $f_{|W}$ lie close to each 
other. This is not a problem here, but in general if we want to avoid 
this situation we need to impose one additional transversality condition
on $f$. Namely, we must require the uniform transversality to $0$ of 
$\partial(f_{|W})$, which is easily obtained by imitating Donaldson's
argument from \cite{D2}. Another situation in which this property naturally
becomes satisfied is the one described in \S 5.
\end{rem}

Given a braided curve $D\subset\CP^2$ of degree $d$ described by a braid
factorization as in (7), and given a monodromy map $\theta_{n-1}$ as in (9), 
certain compatibility conditions need to hold between them in order to ensure 
the existence of a $\CP^2$-valued map with critical curve $D$ and monodromy
$\theta_{n-1}$. Namely, $\theta_{n-1}$ must factor through
$\pi_1(\C^2-D)$, and the fibration must behave in accordance with the 
expected models near the special points of $D$.
We introduce the following definition summarizing these compatibility
properties~:

\begin{defn} A geometric $(n-1)$-dimensional monodromy representation 
associated to a braided curve $D\subset\CP^2$ is a surjective group 
homomorphism $\theta_{n-1}$ from the free group 
$\pi_1(\C-\{q_1,\dots,q_d\})=F_d$ to a symplectic mapping class group
$\mathrm{Map}^\omega(\Sigma^{2n-4},Z^{2n-6})$,
mapping the geometric generators $\gamma_i$ (and thus also the $\gamma_i*Q_j$)
to positive Dehn twists and such that\medskip

$\theta_{n-1}(\gamma_1\dots\gamma_d)=\delta_Z,$

$\theta_{n-1}(\gamma_1*Q_j)=\theta_{n-1}(\gamma_2*Q_j)$ if $r_j=1$,

$\theta_{n-1}(\gamma_1*Q_j)$ and $\theta_{n-1}(\gamma_2*Q_j)$ are twists 
along disjoint Lagrangian spheres if $r_j=\pm 2$,

$\theta_{n-1}(\gamma_1*Q_j)$ and $\theta_{n-1}(\gamma_2*Q_j)$ are twists along
Lagrangian spheres transversely intersecting in one point if $r_j=3$.
\end{defn}

As in the four-dimensional case, $\theta_{n-1}$ remains unchanged and the 
compatibility conditions are preserved when the braid factorization 
defining $D$ is affected by a Hurwitz move. However, when all factors in
the braid factorization are simultaneously conjugated
by a certain braid $Q\in B_d$, the system of geometric generators
$\gamma_1,\dots,\gamma_d$ changes accordingly, and so the geometric monodromy 
representation $\theta_{n-1}$ should be replaced by $\theta_{n-1}\circ Q_*$, 
where $Q_*$ is the automorphism of $F_d$ induced by the braid $Q$. For
example, conjugating the braid factorization by one of the generating
half-twists in $B_d$ affects the monodromy $\theta_{n-1}$ of the restricted
pencil by a Hurwitz move.
\medskip

One easily checks that, given a symplectic braided curve $D\subset\CP^2$ 
and a compatible monodromy representation $\theta_{n-1}:F_d\to
\mathrm{Map}^\omega(\Sigma,Z)$, it is possible
to recover a compact $2n$-manifold $X$ and a map
$f:X-Z\to\CP^2$ in a canonical way up to smooth isotopy. Moreover, 
it is actually possible to endow $X$ with a symplectic structure, 
canonically up to symplectic isotopy. Indeed, by first applying Theorem 2.2
to the monodromy map $\theta_{n-1}$ we can recover a canonical symplectic 
structure on the total space $W$ of the restricted Lefschetz pencil~; 
furthermore, as will be shown in \S 4.4 below,
the braid monodromy of $D$ and the compatible monodromy representation 
$\theta_{n-1}$ determine on $X$ a structure of Lefschetz pencil with generic
fiber $W$ and base set $\Sigma$, which implies by a second application of
Theorem 2.2 that $X$ carries a canonical symplectic structure. The same
result can also be obtained more directly, by adapting the statement and
proof of Theorem 2.2 to the case of $\CP^2$-valued maps.
\medskip

As in the four-dimensional case, we can naturally define symplectic
invariants arising from the quasiholomorphic maps constructed in
Theorem 4.1. However, we again need to take into account the
possible presence of negative self-intersections in the critical 
curves of these maps.
Therefore, the braid factorizations we obtain 
are only canonical up to global conjugation, Hurwitz equivalence, and
pair cancellations or creations. As in the four-dimensional case, a
pair creation operation (inserting two mutually inverse factors anywhere
in the braid factorization) is only allowed if the new factorization 
remains compatible with the monodromy representation $\theta_{n-1}$, 
i.e.\ if $\theta_{n-1}$ maps the two corresponding geometric generators
to Dehn twists along disjoint Lagrangian spheres.

With this understood, we can introduce a notion of m-equivalence as in
Definition 3.5. The following result then holds~:

\begin{thm}
The braid factorizations and geometric monodromy 
representations associated to the quasiholomorphic maps to $\CP^2$ 
obtained in Theorem 4.1 are, for $k\gg 0$, canonical up to m-equivalence, 
and define symplectic invariants of $(X^{2n},\omega)$.

Conversely, the data consisting of a braid factorization and a geometric
$(n-1)$-dimensional monodromy representation, or a m-equivalence class of 
such data, determines a symplectic $2n$-manifold in a canonical way up to 
symplectomorphism.
\end{thm}

\begin{rem}
The invariants studied in this section are a very natural generalization
of those defined in \S 3.2 for 4-manifolds. Namely, when $\dim X=4$, we
naturally get that $Z=\emptyset$ and $\dim \Sigma=0$, i.e.\ the
generic fiber $\Sigma$ consists of a finite number of points, as
expected for a branched covering map. In particular, the mapping class 
group $\mathrm{Map}(\Sigma)$ of the 0-manifold $\Sigma$ is in fact the
symmetric group of order $\mathrm{card}(\Sigma)$. Finally, a Lagrangian
0-sphere in $\Sigma$ is just a pair of points of $\Sigma$, and the
associated Dehn twist is simply the corresponding transposition.
With this correspondence, the results of \S 3 are the exact four-dimensional
counterparts of those described here.
\end{rem}

\subsection{Quasiholomorphic maps and symplectic Lefschetz pencils}

Consider again a symplectic manifold $(X^{2n},\omega)$ and let
$f:X-Z\to\CP^2$ be a map with the same topological properties as those
obtained by Theorem 4.1 from sections of $L^{\otimes k}$ for $k$ large 
enough. As in the four-dimensional case, the $\CP^1$-valued map 
$\pi\circ f$ defines a Lefschetz pencil structure on $X$, obtained by 
lifting via $f$ a pencil of lines on $\CP^2$. The base set of this pencil 
is the fiber of $f$ above the pole $(0\!:\!0\!:\!1)$ of the projection 
$\pi$.

In fact, starting from the quasiholomorphic maps $f_k$ given by Theorem 4.1,
the symplectic Lefschetz pencils $\pi\circ f_k$ coincide for $k\gg 0$
with those obtained by Donaldson in \cite{D2} and described in \S 2~;
calling $s_k^0,s_k^1,s_k^2$ the sections of $L^{\otimes k}$
defining $f_k$, the Lefschetz pencil $\pi\circ f_k$ is the one induced by the 
sections $s_k^0$ and $s_k^1$. 

Therefore, as in the case of a 4-manifold, the invariants described in 
\S 4.3 (braid factorization and $(n-1)$-dimensional geometric monodromy 
representation) completely determine those discussed in \S 2 (factorizations 
in mapping class groups). Once again, the topological description of the 
relation between quasiholomorphic maps and Lefschetz pencils involves a 
subgroup of {\it $\theta_{n-1}$-liftable braids} in the braid group, 
and a group homomorphism from this subgroup to a mapping class group.
\medskip

Consider a symplectic braided curve $D\subset\CP^2$, described by its braid
monodromy $\rho_n:\pi_1(\C-\{p_1,\dots,p_r\})\to B_d$, and a compatible
$(n-1)$-dimensional monodromy representation
$\theta_{n-1}:F_d=\pi_1(\C-\{q_1,\dots,q_d\})\to
\mathrm{Map}^\omega(\Sigma^{2n-4}, Z^{2n-6})$. Then we can make the
following definition~:

\begin{defn} 
The subgroup $B_d^0(\theta_{n-1})$ of liftable braids is the set of all 
braids $Q\in B_d$ such that $\theta_{n-1}\circ Q_*=\theta_{n-1}$, where
$Q_*\in\mathrm{Aut}(F_d)$ is the automorphism induced by the braid $Q$ on
$\pi_1(\C-\{q_1,\dots,q_d\})$.
\end{defn}

A topological definition of $B_d^0(\theta_{n-1})$ can also be given in 
terms of universal fibrations and coverings of configuration spaces, 
similarly to the description in \S 3.3. 

More importantly, denote by $W$ the total space of the symplectic Lefschetz
pencil $LP(\theta_{n-1})$ with generic fiber $\Sigma$ and monodromy 
$\theta_{n-1}$. For example, if $\rho_n$ and $\theta_{n-1}$ are the monodromy
morphisms associated to a quasiholomorphic map given by sections 
$s_k^0,s_k^1,s_k^2$ of $L^{\otimes k}$ over $X$, then $W$ is the smooth 
symplectic hypersurface in $X$ given by the equation $s_k^0=0$~; indeed, as 
seen in \S 4.3, this hypersurface carries a Lefschetz pencil structure 
with generic fiber $\Sigma$, induced by $s_k^1$ and $s_k^2$,  and the 
monodromy of this restricted pencil is precisely $\theta_{n-1}$. 
A braid $Q\in B_d$ can be viewed as a motion of the critical set 
$\{q_1,\dots,q_d\}$ of the Lefschetz pencil $LP(\theta_{n-1})$~; after
this motion we obtain a new Lefschetz pencil with monodromy
$\theta_{n-1}\circ Q_*$. So the subgroup $B_d^0(\theta_{n-1})$ precisely 
consists of those braids which preserve the monodromy of the Lefschetz
pencil $LP(\theta_{n-1})$.

Viewing braids as compactly supported symplectomorphisms of 
the plane preserving $\{q_1,\dots,q_d\}$, the fact that $Q$ belongs to
$B_d^0(\theta_{n-1})$ means that it can be lifted via the Lefschetz pencil
map $W-Z\to\CP^1$ to a symplectomorphism of $W$. Since the monodromy of
the pencil $LP(\theta_{n-1})$ preserves a neighborhood of the base set $Z$,
the lift to $W$ of the braid $Q$ coincides with the identity over a 
neighborhood of $Z$. Even better, because $Q$ is compactly supported, its lift
to $W$ coincides with $\mathrm{Id}$ near the fiber above the
point at infinity in $\CP^1$, which can be identified with $\Sigma$.
Therefore, the lift of $Q$ to $W$ is a well-defined element of the mapping
class group $\mathrm{Map}^\omega(W,\Sigma)$, which we call 
$(\theta_{n-1})_*(Q)$. This construction defines a group homomorphism
$$(\theta_{n-1})_*:B_d^0(\theta_{n-1})\to\mathrm{Map}^\omega(W^{2n-2},
\Sigma^{2n-4}).$$

Since the geometric monodromy representation $\theta_{n-1}$ is compatible
with the braided curve $D\subset\CP^2$, the image of the braid 
monodromy homomorphism $\rho_n:\pi_1(\C-\{p_1,\dots,p_r\})\to B_d$ describing
$D$ is entirely contained in $B^0_d(\theta_{n-1})$. Indeed, it follows from
Definition 4.1 that $\theta_{n-1}$ factors through $\pi_1(\C^2-D)$, on which
the braids of $\mathrm{Im}\,\rho_n$ act trivially.
As a consequence, we can use the group homomorphism $(\theta_{n-1})_*$ in
order to obtain, from the braid monodromy $\rho_n$, a group homomorphism
$$\theta_n=(\theta_{n-1})_*\circ\rho_n:\pi_1(\C-\{p_1,\dots,p_r\})\to
\mathrm{Map}^\omega(W,\Sigma).$$
If $\rho_n$ and $\theta_{n-1}$ describe the monodromy of a $\CP^2$-valued
map $f$, then $\theta_n$ is by construction the monodromy of the 
corresponding Lefschetz pencil $\pi\circ f$. Therefore, the following 
result holds~:

\begin{prop}
Let $f:X-Z\to\CP^2$ be one of the quasiholomorphic maps of Theorem 4.1.
Let $D\subset\CP^2$ be its critical curve, and denote by
$\rho_n:\pi_1(\C-\{p_1,\dots,p_r\})\to B_d^0(\theta_{n-1})$ and
$\theta:F_d\to\mathrm{Map}^\omega(\Sigma,Z)$ be the corresponding monodromies.
Then the monodromy map $\theta_n:\pi_1(\C-\{p_1,\dots,p_r\})\to 
\mathrm{Map}^\omega(W,\Sigma)$ of the Lefschetz pencil $\pi\circ f$ 
is given by the identity $\theta_n=(\theta_{n-1})_*\circ\rho_n$.

In particular, for $k\gg 0$ the symplectic invariants given by
Theorem 2.1 are obtained in this manner from those defined in Theorem 4.3.
\end{prop}

As in the four-dimensional case, all the factors of degree $\pm 2$ or $3$
in the braid monodromy (corresponding to the cusps and nodes of $D$) lie in 
the kernel of $(\theta_{n-1})_*$~; the only terms which
contribute non-trivially to the pencil monodromy $\theta_n$
are those arising from the tangency points of the branch curve $D$, and each
of these contributions is a Dehn twist. 

More precisely, the image in $\mathrm{Map}^\omega(W,\Sigma)$
of a half-twist $Q\in B^0_d(\theta_{n-1})$ arising as the braid monodromy
around a tangency point of $D$ can be constructed as follows. 
Consider the Lefschetz pencil $LP(\theta_{n-1})$ with total space $W$, 
generic fiber $\Sigma$, critical levels $q_1,\dots,q_d$ and monodromy 
$\theta_{n-1}$. Call $\gamma$ the path joining two of the points 
$q_1,\dots,q_d$ (e.g., $q_{i_1}$ and $q_{i_2}$) and naturally associated to 
the half-twist $Q$ (the path along which the twisting occurs). By 
Definition 4.1, the monodromies of $LP(\theta_{n-1})$ around the two end points
$q_{i_1}$ and $q_{i_2}$ are the same Dehn twists (using $\gamma$ to identify 
the two singular fibers). Even better, in this context one easily shows that
the vanishing cycles at the two end points of $\gamma$ are isotopic 
Lagrangian spheres in $\Sigma$. Then it follows from the work of Donaldson
and Seidel that, above the path $\gamma$, one can find a
Lagrangian sphere $L=S^{n-1}\subset W$, joining the singular points of
the fibers above $q_{i_1}$ and $q_{i_2}$, and intersecting each fiber
inbetween in a Lagrangian sphere $S^{n-2}$. The element $(\theta_{n-1})_*(Q)$ 
in $\mathrm{Map}^\omega(W,\Sigma)$ is the positive Dehn twist along
the Lagrangian sphere $L$.

\begin{rem}
Let $(X^{2n},\omega)$ be a compact symplectic manifold, and consider the
symplectic Lefschetz pencils given by Donaldson's result (Theorem 2.1) from
pairs of sections of $L^{\otimes k}$ for $k\gg 0$~; the monodromy of these 
Lefschetz pencils consists of generalized Dehn twists around Lagrangian 
$(n-1)$-spheres in the generic fiber $W_k$. It follows from 
Proposition 4.4 that these Lagrangian spheres are not arbitrary. 
Indeed, they can all be obtained by endowing $W_k$ with a structure of 
symplectic Lefschetz pencil induced by two sections of $L^{\otimes k}$ 
(the existence of such a structure follows from the results of this
section), and by looking for Lagrangian $(n-1)$-spheres which join two 
mutually isotopic vanishing cycles of this pencil above a path in the base.

As observed by Seidel, this remarkable structure of vanishing cycles makes
it possible to hope for a purely combinatorial description of Lagrangian
Floer homology~: one can try to use the structure of vanishing 
cycles in a $2n$-dimensional Lefschetz pencil to reduce things first to the
$2n-2$-dimensional case, and then by induction eventually to the case of
$0$-manifolds, in which the calculations are purely combinatorial.
\end{rem}

\section{Complete linear systems and dimensional induction}

We now show how the results of \S 4 can be used in order to reduce in
principle the classification of compact symplectic manifolds to a purely
combinatorial problem.

\begin{defn}
Let $(X^{2n},\omega)$ be a compact symplectic manifold. We say that 
asymptotically holomorphic $(n+1)$-tuples of sections of $L^{\otimes k}$
define {\it 3-good complete linear systems} on $X$ if, for large values
of $k$, these sections $s_0,\dots,s_n\in\Gamma(L^{\otimes k})$ satisfy the 
following properties~:

$(a)$ for $0\le r\le n-1$, the section $(s_{r+1},\dots,s_n)$ of
$\C^{n-r}\otimes L^{\otimes k}$ satisfies a uniform transverslity property,
and its zero set $\Sigma_r=\{s_{r+1}=\dots=s_n=0\}$ is a smooth
symplectic submanifold of dimension $2r$ in $X$. We also define
$\Sigma_n=X$ and $\Sigma_{-1}=\emptyset$~;

$(b)$ for $1\le r\le n$, the pair of sections
$(s_r,s_{r-1})\in\Gamma(\C^2\otimes L^{\otimes k})$ defines a structure of
symplectic Lefschetz pencil on $\Sigma_r$, with generic fiber $\Sigma_{r-1}$
and base set $\Sigma_{r-2}$~;

$(c)$ for $2\le r\le n$, the triple of sections
$(s_r,s_{r-1},s_{r-2})\in\Gamma(\C^3\otimes L^{\otimes k})$ defines a
quasiholomorphic map from $\Sigma_r$ to $\CP^2$, with generic fiber
$\Sigma_{r-2}$ and base set $\Sigma_{r-3}$.
\end{defn}

One can think of a 3-good complete linear system in the following way.
First, the two sections $s_n$ and $s_{n-1}$ define a Lefschetz pencil
structure on $X$. By adding the section $s_{n-2}$, this structure is
refined into a quasiholomorphic map to $\CP^2$. As observed in \S 4,
by restricting to the hypersurface $\Sigma_{n-1}$ we get a symplectic
Lefschetz pencil defined by $s_{n-1}$ and $s_{n-2}$. This structure is in
turn refined into a quasiholomorphic map by adding the section $s_{n-3}$~;
and so on.

Note that, except for the case $r=1$, part $(b)$ of Definition 5.1 is actually
an immediate consequence of part $(c)$, because by composing $\CP^2$-valued 
quasiholomorphic maps with the projection $\pi:\CP^2-\{(0\!:\!0\!:\!1)\}\to
\CP^1$ one always obtains Lefschetz pencils. Also note that, in order to
make sense out of these properties, one implicitly needs to endow the
submanifolds $\Sigma_r$ with $\omega$-compatible almost-complex structures~;
these restricted almost-complex structures can be chosen to differ from the 
almost-complex structure $J$ on $X$ by $O(k^{-1/2})$, so that asymptotic 
holomorphicity and transversality properties are not affected by this
choice.

\begin{thm}
Let $(X^{2n},\omega)$ be a compact symplectic manifold. Then for all large
enough values of $k$ it is possible to find asymptotically holomorphic 
sections of $\C^{n+1}\otimes L^{\otimes k}$ determining 3-good complete
linear systems on $X$. Moreover, for large $k$ these structures are 
canonical up to isotopy and up to cancellations of pairs of nodes in the 
critical curves of the quasiholomorphic $\CP^2$-valued maps.
\end{thm}

\begin{proof} We only give a sketch of the proof of Theorem 5.1.
As usual, we need to obtain two types of properties~: uniform transversality
conditions, which we ensure in the first part of the argument, and 
compatibility conditions, which are obtained by a subsequent perturbation.
As in previous arguments, the various uniform transversality properties are
obtained successively, using the fact that, because transversality is an open
condition, it is preserved by any sufficiently small subsequent perturbations.

The first transversality properties to be obtained are those appearing in 
part $(a)$ of Definition 5.1, i.e.\ the transversality to $0$ of 
$(s_{r+1},\dots,s_n)$ for all $0\le r\le n-1$~; this easy case is 
e.g.\ covered by the main result of \cite{A1}.

One next turns to the transversality conditions arising from the requirement
that the three sections $(s_n,s_{n-1},s_{n-2})$ define quasiholomorphic
maps from $X$ to $\CP^2$~: it follows immediately from the proof of Theorem
4.1 that these properties can be obtained by suitable small perturbations.

Next, we try to modify $s_{n-1}$, $s_{n-2}$ and $s_{n-3}$ in order to ensure
that the restrictions to $\Sigma_{n-1}=s_n^{-1}(0)$ of these three sections
satisfy the transversality properties of Definition 3.2. A general strategy
to handle this kind of situation is to use the following remark (Lemma 6 of
\cite{A2})~: if $\phi$ is a section of a vector bundle $\mathcal{F}$ over
$X$, satisfying a uniform transversality property, and if $W=\phi^{-1}(0)$,
then the uniform transversality to $0$ over $W$ of a section $\xi$ of a 
vector bundle $\mathcal{E}$ is equivalent to the uniform transversality to
$0$ over $X$ of the section $\xi\oplus\phi$ of
$\mathcal{E}\oplus\mathcal{F}$, up to a change in transversality estimates.
This makes it possible to replace all transversality properties to be 
satisfied over submanifolds of $X$ by transversality properties to be
satisfied over $X$ itself~; each property can then be ensured by the
standard type of argument, using the globalization principle to combine
suitably chosen local perturbations (see \cite{Atrans} for more details).

However, in our case the situation is significantly simplified by the fact 
that, no matter how we perturb the sections $s_{n-1}$, $s_{n-2}$ and
$s_{n-3}$, the submanifold $\Sigma_{n-1}$ itself is not affected. Moreover,
the geometry of $\Sigma_{n-1}$ is controlled by the transversality
properties obtained on $s_n$~; for example, a suitable choice
of the constant $\rho>0$ (independent of $k$) ensures that the intersection
of $\Sigma_{n-1}$ with any ball of $g_k$-radius $\rho$ centered at one of
its points is topologically a ball (see e.g.\ Lemma 4 of \cite{A1}).
Therefore, we can actually imitate all steps of the argument used to prove
Theorem 4.1, working with sections of $L^{\otimes k}$ over $\Sigma_{n-1}$.
The localized reference sections of $L^{\otimes k}$ over $\Sigma_{n-1}$ that
we use in the arguments are now chosen to be the restrictions to 
$\Sigma_{n-1}$ of the localized sections $s_{k,x}^\mathrm{ref}$ of
$L^{\otimes k}$ over $X$~; similarly, the approximately holomorphic
local coordinates over $\Sigma_{n-1}$ in which we work are obtained as the
restrictions to $\Sigma_{n-1}$ of local coordinate functions on $X$.
With these two differences understood, we can still construct localized 
perturbations by the same algorithms as in \S 4.1 and, using the standard
globalization argument, achieve the desired transversality properties over 
$\Sigma_{n-1}$. Moreover, all these local perturbations are obtained as
products of the localized reference sections by polynomial functions of the
local coordinates. Therefore, they naturally arise as restrictions to
$\Sigma_{n-1}$ of localized sections of $L^{\otimes k}$ over $X$, and so
we actually obtain well-defined perturbations of the sections $s_{n-1}$,
$s_{n-2}$ and $s_{n-3}$ over $X$ which yield the desired transversality
properties over $\Sigma_{n-1}$.

We can continue similarly by induction on the dimension, until we obtain
the transversality properties required of $s_2$, $s_1$ and $s_0$ over
$\Sigma_2$, and finally the transversality properties required of $s_1$ and
$s_0$ over $\Sigma_1$. Observe that, even though the perturbations performed
over each $\Sigma_r$ result in modifications of the submanifolds $\Sigma_j$ 
($j<r$) lying inside them, these perturbations preserve the
transversality properties of $(s_{j+1},\dots,s_n)$,
and so the submanifolds $\Sigma_j$ retain their smoothness and
symplecticity properties.

We now turn to the second part of the argument, i.e.\ obtaining the desired
compatibility conditions. First observe that the proof of Theorem 4.1
shows how, by a perturbation of $s_n$, $s_{n-1}$ and $s_{n-2}$ smaller than
$O(k^{-1/2})$, we can ensure that the various compatibility properties of
Definition 3.2 are satisfied by the $\CP^2$-valued map $f_n$ defined by 
these three sections.

Next, we proceed to perturb $f_{n-1}=(s_{n-1}\!:\!s_{n-2}\!:\!s_{n-3})$ over a 
neighborhood of its ramification curve $R_{n-1}\subset\Sigma_{n-1}$, in order 
to obtain the required compatibility properties for $f_{n-1}$, but without 
losing those previously achieved for $f_n$ near its ramification curve
$R_n\subset X$. For this purpose, we first show that the curve $R_n$
satisfies a uniform transversality property with respect to
the hypersurface $\Sigma_{n-1}$ in $X$. 

The only way in which $R_n$
can fail to be uniformly transverse to $\Sigma_{n-1}$ is if 
$\partial(\pi\circ f_{n|R_n})$ becomes small at a point of $R_n$ near 
$\Sigma_{n-1}$. Because $f_n$ satisfies property $(6)$ in Definition 3.2, 
this can only happen if a cusp point or a tangency point of $f_n$ lies 
close to $\Sigma_{n-1}$. However, property $(7)$ of Definition 3.2 implies
that this point cannot belong to $\Sigma_{n-1}$. Therefore, two of the
intersection points of $R_n$ with $\Sigma_{n-1}$ must lie close to each 
other. Observe that the points of
$R_n\cap\Sigma_{n-1}$ are precisely the critical points of the Lefschetz
pencil induced on $\Sigma_{n-1}$ by $s_{n-1}$ and $s_{n-2}$, i.e.\ the
tangency points of the map $f_{n-1}$. The
transversality properties already obtained for $f_{n-1}$ imply that two 
tangency points cannot lie close to each other~; we get a contradiction,
so the cusps and tangencies of $f_n$ must lie far away from $\Sigma_{n-1}$, 
and $R_n$ and $\Sigma_{n-1}$ are mutually transverse.

This implies in particular that a small perturbation of $s_{n-1}$, $s_{n-2}$ 
and $s_{n-3}$ localized near $\Sigma_{n-1}$ cannot affect properties $(4')$
and $(6')$ for $f_n$, and also that the only place where perturbing
$f_{n-1}$ might affect $f_n$ is near the tangency points of $f_{n-1}$.

We now consider the set
$\mathcal{C}_{n-1}\cup\mathcal{T}_{n-1}\cup\mathcal{I}_{n-1}$ of points
where we need to ensure properties $(4')$, $(6')$ and $(8')$ for $f_{n-1}$.
The first step is as usual to perturb $J$ into an almost-complex structure
which is integrable near these points~; once this is done, we perturb
$f_{n-1}$ to make it locally holomorphic with respect to this almost-complex
structure. 

We start by considering a point $x\in\mathcal{C}_{n-1}\cup\mathcal{I}_{n-1}$, 
where the issue of preserving properties of $f_n$ does not arise. We follow 
the argument in \S 4.1 of \cite{A2}. First, it is possible to perturb the
almost-complex structure $J$ over a neighborhood of $x$ in $X$ in order to
obtain an almost-complex structure $\tilde{J}$ which differs from $J$ by
$O(k^{-1/2})$ and is integrable over a small ball centered at $x$.
Recall from \cite{A2} that $\tilde{J}$ is obtained by choosing approximately
holomorphic coordinates on $X$ and using them to pull back the standard 
complex structure of $\C^n$~; a cut-off function is used to splice $J$ with
this locally defined integrable structure. Since we can choose the local 
coordinates in such a way that a local equation of $\Sigma_{n-1}$ is $z_n=0$, 
we can easily ensure that $\Sigma_{n-1}$ is, over a small neighborhood of
$x$, a $\tilde{J}$-holomorphic submanifold of $X$. Next, we can 
perturb the sections $s_{n-1},s_{n-2},s_{n-3}$ of $L^{\otimes k}$ by 
$O(k^{-1/2})$ in order to make the projective map defined by them
$\tilde{J}$-holomorphic over a neighborhood of $x$ in $X$ (see \cite{A2}). 
This holomorphicity property remains true for the restrictions to the locally
$\tilde{J}$-holomorphic submanifold $\Sigma_{n-1}$. So, we have obtained
the desired compatibility property near $x$.

We now consider the case of a point $x\in\mathcal{T}_{n-1}$, where we need
to obtain property $(6')$ for $f_{n-1}$ while preserving property $(8')$ for
$f_n$. We first observe that, by the construction of the previous step
(getting property $(8')$ for $f_n$ at $x$), we have a readily available 
almost-complex structure $\tilde{J}$ integrable over a neighborhood of $x$ 
in $X$. In particular, by construction $f_n$ is locally 
$\tilde{J}$-holomorphic and $\Sigma_{n-1}$ is locally a 
$\tilde{J}$-holomorphic submanifold of $X$. We next try to make the
projective map $f_{n-1}$ holomorphic over a neighborhood of $x$, using once
again the argument of \cite{A2}. The key observation here is that, because
one of the sections $s_{n-1}$ and $s_{n-2}$ is bounded from below at $x$,
we can reduce to a $\C^2$-valued map whose first component is already 
holomorphic. Therefore, the perturbation process described in \cite{A2} only
affects $s_{n-3}$, while the two other sections are preserved. This means
that we can ensure the local $\tilde{J}$-holomorphicity of $f_{n-1}$ without
affecting $f_n$.

It is easy to combine the various localized perturbations performed near
each point of $\mathcal{C}_{n-1}\cup\mathcal{T}_{n-1}\cup\mathcal{I}_{n-1}$~;
this yields properties $(4')$, $(6')$ and $(8')$ of Definition 3.2
for $f_{n-1}$.

We now use a generically chosen small perturbation of $s_{n-1}$, $s_{n-2}$
and $s_{n-3}$ in order to ensure property $(7)$, i.e.\ the
self-transversality of the critical curve of $f_{n-1}$. It is important to
observe that, because $f_n$ satisfies property $(7)$, the images by the 
projective map $(s_{n-1}\!:\!s_{n-2})$ of the points of
$R_n\cap\Sigma_{n-1}=\mathcal{I}_n=\mathcal{T}_{n-1}$ are all distinct from
each other, and because $f_n$ satisfies property $(5)$ they are also distinct 
from $(0\!:\!1)$. Therefore, we can choose a perturbation which vanishes
identically over a neighborhood of $\mathcal{T}_{n-1}$~; this makes it
possible to obtain property $(7)$ for $f_{n-1}$ without losing any
property of $f_n$.

Finally, by the process described in \S 4.2 of \cite{A2} we construct a 
perturbation yielding property $(3')$ along the critical curve of $f_{n-1}$~;
this perturbation is originally defined only for the restrictions to
$\Sigma_{n-1}$ but it can easily be extended outside of $\Sigma_{n-1}$ by 
using a cut-off function. The two important properties of this perturbation
are the following~: first, it vanishes identically near the points where 
$f_{n-1}$ has already been made $\tilde{J}$-holomorphic, and in particular 
near the points of $\mathcal{T}_{n-1}$~; therefore, none of the properties 
of $f_n$ are affected, and properties $(4')$, $(6')$ and $(8')$ of $f_{n-1}$
are not affected either. Secondly, this perturbation does not modify
the critical curve of $f_{n-1}$ nor its image, so property $(7)$ is
preserved. We have therefore obtained all desired properties for $f_{n-1}$.

We can continue similarly by induction on the dimension, until all required
compatibility properties are satisfied. Observe that, because the 
ramification curve of $f_r$ remains
away from its fiber at infinity $\Sigma_{r-2}$, we do not need
to worry about the possible effects on $f_r$ of perturbations of
$f_{r-2}$. Therefore, the argument remains the same at each step, and we
can complete the proof of the existence statement in Theorem 5.1 in this way.

The proof of the uniqueness statement relies, as usual, on the extension of
the whole construction to one-parameter families~; this is easily done by
following the same ideas as in previous arguments.
\end{proof}

The structures of 3-good complete linear systems given by Theorem 5.1 are
extremely rich, and lead to interesting invariants of compact symplectic
manifolds. Indeed, recall from Definition 5.1 that, for $1\le r\le n$, the
sections $s_r$ and $s_{r-1}$ define a symplectic Lefschetz pencil structure
on $\Sigma_r$, with generic fiber $\Sigma_{r-1}$ and base set $\Sigma_{r-2}$.
The monodromy of this pencil is given by a group homomorphism
\begin{equation}\theta_r:\pi_1(\C-\{p_1,\dots,p_{d_r}\})\to
\mathrm{Map}^\omega(\Sigma_{r-1},\Sigma_{r-2}).\end{equation}

Moreover, for $2\le r\le n$, the sections $s_r$, $s_{r-1}$ and $s_{r-2}$ 
define a quasiholomorphic map from $\Sigma_r-\Sigma_{r-3}$ to $\CP^2$, 
with generic fiber $\Sigma_{r-2}$. Denote by $D_r\subset\CP^2$ the critical
curve of this map, and let $d_{r-1}=\deg D_r$. As shown in \S 4.3, we obtain
two monodromy morphisms~: on one hand, the braid monodromy homomorphism
characterizing $D_r$,
\begin{equation}\rho_r:\pi_1(\C-\{p_1,\dots,p_{s_r}\})\to B_{d_{r-1}},
\end{equation}
and on the other hand, a compatible $(r-1)$-dimensional monodromy
representation, which was shown in \S 4.3 to be none other than
$$\theta_{r-1}:\pi_1(\C-\{p_1,\dots,p_{d_{r-1}}\})\to
\mathrm{Map}^\omega(\Sigma_{r-2},\Sigma_{r-3}).$$

Finally, it was shown in \S 4.4 that $\mathrm{Im}(\rho_r)\subseteq 
B_{d_{r-1}}^0(\theta_{r-1})$, and that the various monodromies are related
to each other by the identity \begin{equation}
\theta_r=(\theta_{r-1})_*\circ\rho_r.\end{equation}

In particular, the manifold $X$ is completely characterized by the braid
monodromies $\rho_2,\dots,\rho_n$ and by the map
$\theta_1$ with values in $\mathrm{Map}^\omega(\Sigma_0,\emptyset)$, which
is a symmetric group~; this data is sufficient to successively
reconstruct all morphisms $\theta_r$ and all submanifolds $\Sigma_r$ by
inductively using equation $(12)$. 

In other words, a symplectic $2n$-manifold is characterized by $n-2$ braid 
factorizations and a word in a symmetric group~; or, stopping at
$\theta_2$, we can also consider $n-3$ braid factorizations and a
word in the mapping class group of a Riemann surface.

These results can be summarized by the following theorem~:

\begin{thm}
The braid monodromies $\rho_2,\dots,\rho_n$ and the symmetric group
representation $\theta_1$ associated to the 3-good complete linear systems
obtained in Theorem 5.1 are, for $k\gg 0$, canonical up to m-equivalence, 
and define symplectic invariants of $(X^{2n},\omega)$.

Conversely, the data consisting of several braid factorizations and a
symmetric group representation satisfying suitable compatibility
conditions, or a m-equivalence class of 
such data, determines a symplectic $2n$-manifold in a canonical way up to 
symplectomorphism.
\end{thm}

In principle, this result reduces the study of compact symplectic manifolds 
to purely combinatorial questions about braid groups and symmetric groups~; 
however, the invariants it introduces are probably quite difficult to compute
as soon as one considers examples which are not complex algebraic.
Nevertheless, it seems that this construction should be very helpful in
improving our understanding of the topology of Lefschetz pencils in 
dimensions greater than $4$.

\end{document}